\documentclass[10pt]{amsart}
\usepackage{amsmath,amscd,amssymb,verbatim,paralist}
\usepackage[percent]{overpic}
\graphicspath{{figs/}}

\setdefaultleftmargin{2.5em}{2.2em}{1.87em}{1.7em}{1em}{1em}
\newcommand{\ZZ}{\mathbb{Z}}
\newcommand{\RR}{\mathbb{R}}
\newcommand{\CC}{\mathbb{C}}

\newcommand{\AAA}{{\mathcal{A}}}
\newcommand{\CCC}{{\mathcal{C}}}
\newcommand{\RRR}{{\mathcal{R}}}
\newcommand{\MMM}{{\mathcal{M}}}

\newcommand{\sym}{\mathfrak{S}}

\newcommand{\Cov}{\operatorname{Cov}}
\newcommand{\rank}{\operatorname{rank}}
\newcommand{\conn}{\operatorname{conn}}
\newcommand{\sgn}{\operatorname{sgn}}

\newcommand{\symmdiff}{\bigtriangleup}

\newtheorem{theorem}{Theorem}[section]
\newtheorem{corollary}[theorem]{Corollary}
\newtheorem{proposition}[theorem]{Proposition}
\newtheorem{conjecture}[theorem]{Conjecture}

\newtheoremstyle{defn}{1.2ex}{1.2ex}{}{}{}{.}{.5em}%
{\textbf{\thmname{#1}\thmnumber{ #2}}\thmnote{\emph{ (#3)}}}
\theoremstyle{defn}
\newtheorem{defn}[theorem]{Definition}

\newtheorem{example}[theorem]{Example}
\newtheorem{remark}[theorem]{Remark}
\newtheorem{question}[theorem]{Question}

\numberwithin{figure}{section}
\numberwithin{table}{section}

\begin{document}
\title{Diameter of graphs of reduced words and galleries}

\author{Victor Reiner}
\address{School of Mathematics\\
University of Minnesota\\
Minneapolis, MN 55455\\
USA}
\email{reiner@math.umn.edu}

\author{Yuval Roichman}
\address{Department of Mathematics\\
Bar-Ilan University\\
52900 Ramat-Gan\\
Israel}
\email{yuvalr@math.biu.ac.il}


\keywords{Coxeter group, reduced words, supersolvable, hyperplane arrangement,
weak order, reflection order, cellular string, zonotope, 
monotone path, diameter}

\subjclass[2000]{20F55,20F05}

\thanks{
First author supported by NSF grant DMS--0245379. Second author
supported in part by the Israel Science Foundation grant \# 947/04.
}

\begin{abstract}
For finite reflection groups of types $A$ and $B$, we
determine the diameter of the graph whose vertices are reduced words
for the longest element and whose edges are braid relations.
This is deduced from a more general theorem that applies to
supersolvable hyperplane arrangements.
%
\end{abstract}

\maketitle


\section{Introduction}
\label{section:intro}

The symmetric group $W=\sym_n$ on $n$ letters has a well-known
Coxeter presentation, with generating set $S=\{s_1,\ldots,s_{n-1}\}$
consisting of the adjacent transpositions $s_i=(i,i+1)$, satisfying
the \emph{braid relations}
$$
\begin{array}{ccl}
\text{(i)}  & s_i s_j = s_j s_i & \text{ for }|i - j| \geq 2,\\
\text{(ii)} & s_i s_{i+1} s_i = s_{i+1} s_i s_{i+1} & \text{ for }1 \leq i \leq n-2,
\end{array}
$$
together with the condition that
each $s_i$ is an involution.  Given any $w$ in $W$, 
a \emph{reduced decomposition} for $w$
is a sequence $(s_{i_1},\ldots,s_{i_\ell})$ of the generators $S$
for which $w= s_{i_1} \cdots s_{i_\ell}$ that attains the minimum possible
length $\ell=:\ell(w)$.

\begin{figure}[h]
\centerline{\begin{overpic}[width=0.9\textwidth]{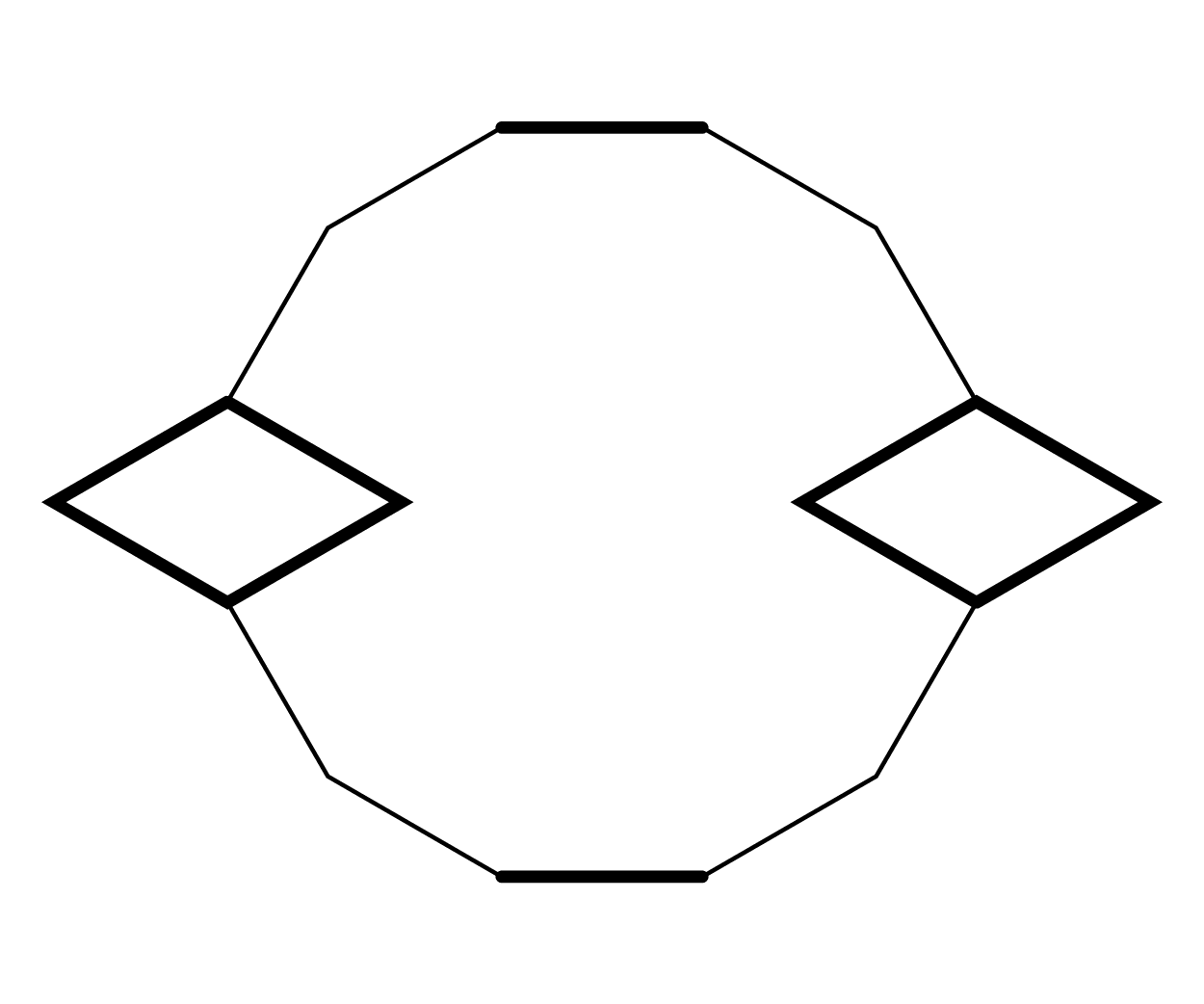}
\put(31,7){$121321$} \put(59,7){$123121$}
\put(17,16){$212321$} \put(73,16){$123212$}
\put(9,30){$213231$} \put(82,30){$132312$}
\put(-7,41){$231231$} \put(97,41){$132132$}
\put(35,41){$213213$} \put(55,41){$312312$}
\put(9,51){$231213$} \put(82,51){$312132$}
\put(17,65){$232123$} \put(73,65){$321232$}
\put(31,74){$323123$} \put(59,74){$321323$}
\end{overpic}}
\caption{The graph $G(w_0)$ for $W=\sym_4$.}\label{fig:S4-example-figure}
\end{figure}

There is a well-studied graph $G(w)$ whose vertex set is
the set $\RRR(w)$ of all reduced decompositions of $w$,
and whose edges correspond to the applicable braid relations (i) and (ii) above.
A theorem of Tits~\cite{Tits-words} (see also \cite[Theorem 3.3.(ii)]{BjornerBrenti})
says that for any finite Coxeter group $(W,S)$
and any $w$ in $W$, this graph $G(w)$ is {\it connected}.
A particularly interesting special case occurs when $w$ is the unique longest
element $w_0$ of $W$.  For $W=\sym_4$, the graph $G(w_0)$ is illustrated
in Figure~\ref{fig:S4-example-figure}, where each reduced word
is abbreviated by its subscript sequence (e.g., $121321$ for
$(s_1,s_2,s_1,s_3,s_2,s_1)$), and with braid relations of type (i) 
darkened.

The graph $G(w_0)$ and some of its generalizations were shown
to have further {\it graph-theoretic} connectivity in work of
Athanasiadis, Edelman and Reiner~\cite{AthanasiadisEdelmanReiner},
and Athanasiadis and Santos~\cite{AthanasiadisSantos}.  This was
motivated by earlier {\it topological} connectivity results surrounding
a closely related poset, appearing first in
a conjecture of Baues~\cite{Baues} on loop spaces, which was proven in
work of Billera, Kapranov, and Sturmfels~\cite{BilleraKapranovSturmfels}
and Bj\"orner~\cite{Bjorner}.
We also mention here a few ancillary
results about the graph $G(w_0)$.
Tits~\cite{Tits-local} gave explicit generators for its fundamental group.
Stanley~\cite{Stanley-words} was the first to show that its vertex set $\RRR(w_0)$ is equinumerous with the standard Young tableaux of
shape $(n-1,n-2,\ldots,2,1)$.
In~\cite{Reiner}, the average degree of a vertex of $G(w_0)$
with respect to only the edges of type (ii) was shown to be $1$.
Manin and Schechtman~\cite{ManinSchechtman}, Ziegler~\cite{Ziegler},
Felsner~\cite{Felsner}, and Shapiro--Shapiro--Vainshtein~\cite{ShapiroShapiroVainshtein}
have studied, in the guise of the \emph{higher Bruhat order} $B(n,2)$,
the quotient graph of $G(w_0)$ in which one contracts down all its edges of type (i).

However, the \emph{diameter} of $G(w_0)$, seems to have been considered only very recently.
Autord and Dehornoy \cite[Proposition 1]{AutordDehornoy}
show that for $W=\sym_n$, the diameter of $G(w_0)$
grows asymptotically in $n$ as a constant times $n^4$.
Our main result, Theorem~\ref{thm:supersolvable-theorem},
shows the diameter is exactly
$
\frac{1}{24}(n-2) (n-1) n (3n-5),
$
which is the number of codimension-two subspaces arising
as intersections of two hyperplanes $x_i=x_j$ in the
{\it reflection arrangement} associated to $W=\sym_n$.

\subsection{Diameter of $G(w_0)$ and supersolvable hyperplane arrangements.}

The graph $G(w_0)$ for any finite Coxeter group $(W,S)$
has a natural generalization to the context of \emph{real
hyperplane arrangements} that first arose in work of Deligne~\cite{Deligne},
and later Salvetti~\cite{Salvetti}, on the topology of the complexified
complements of these arrangements.  This was generalized further
to the context of \emph{oriented matroids}
by Cordovil and Moreira~\cite{CordovilMoreira};
to decrease technicalities and enhance readability, we will
mainly adhere to the language of hyperplane arrangements in this paper.
We review the arrangement viewpoint here in order to state
Theorem~\ref{thm:supersolvable-theorem};
see \cite[\S 4.4, pp.~184--186]{OMbook} and Remark~\ref{rmk:connectivity-credit}
below for further discussion.

Let $\AAA$ be an arrangement of finitely many linear hyperplanes
in $\RR^d$ that is \emph{central} and \emph{essential}, meaning that
$\bigcap_{H \in \AAA} H = \{0\}$.  Let $L=\bigsqcup_{i=0}^d L_i$ 
be its graded poset of intersection subspaces, 
ordered via reverse inclusion.

Define a graph structure $G_1$ on the set
$\CCC$ of chambers of $\AAA$, in which two chambers $c, c'$ are
connected by an edge when they are \emph{separated} by exactly one
hyperplane $H$ in $L_1$.  It is well-known, and will be recalled
in Section~\ref{section:distances}, why this graph $G_1$ always has diameter exactly
$|L_1|=|\AAA|$, that is, the number of hyperplanes.

Now choose a particular base chamber $c_0$, and
let $\RRR$ denote the set of all minimal \emph{galleries} $r$
(that is, geodesics in~$G_1$) from $c_0$ to $-c_0$.
There is a graph structure $G_2$ on this set $\RRR$,
in which two galleries $r, r'$ are connected
by an edge when they are separated (in a sense made precise
in Section~\ref{section:first-two-graphs})
by exactly one codimension-two intersection subspace $X$ in $L_2$.
This graph $G_2$ is known to be 
connected (see Remark~\ref{rmk:connectivity-credit}
below), and it will be shown in Section~\ref{section:distances} that its
diameter is always \emph{at least} $|L_2|$, raising the
following 
question.

\vskip.05in
\noindent
{\bf Main question.}
{\it 
For real hyperplane arrangements $\AAA$ and a choice of base
chamber $c_0$, does the graph $G_2$ of minimal galleries from
$c_0$ to $-c_0$ have diameter $|L_2|$?
}
\vskip.05in
\noindent
Remark~\ref{remark:low-dimension} discusses why the answer to this question is 
{\it affirmative} for  hyperplane arrangements in dimension at most $3$, and 
even for oriented matroids of rank at most $3$, by a result of Cordovil, but
{\it negative} for oriented matroids in rank $4$, by an example of Richter-Gebert.  
After acceptance of this paper for publication, computations by Rob Edman
provided a negative answer for the $4$-dimensional {\it cyclic arrangement} 
of $8$ hyperplanes, that is, the hyperplanes in $\RR^4$ normal to the vectors
$v_i:=(1,t_i,t_i^2,t_i^3)$ for real numbers $t_1 < t_2 < \cdots < t_8$.  Specifically,
in this example, one finds several choices of base chamber $c_0$ (including 
the chamber $c_0$ where vectors have positive dot product with all $8$ of the
$v_i$) for which
the graph $G_2$ of minimal galleries 
from $c_0$ to $-c_0$ has diameter $30>28=\binom{8}{2}=|L_2|$.

Note that when $\AAA$ is the arrangement of reflecting hyperplanes for
a finite real reflection group $W$, the choice of base chamber $c_0$ is
immaterial, as $W$ acts simply transitively on the chambers $\CCC$.
Also, in this case the graph $G_2$ is easily seen to be exactly the graph of reduced
words for the longest element $w_0$ in $W$ described above.

The following main result answers the main question affirmatively
for reflection arrangements of types $A, B$, as well as the more general
\emph{supersolvable arrangements}.  
See Section~\ref{section:supersolvable} for undefined terms
in its statement.

\begin{theorem}
\label{thm:supersolvable-theorem}
When $\AAA$ is a supersolvable hyperplane arrangement,
and the base chamber $c_0$ is chosen incident to a modular flag,
the graph $G_2$ has diameter $|L_2|$.

In particular, for the reflection arrangements of type $A_{n-1}$ and $B_n$
and the dihedral groups $I_2(m)$, the graphs of reduced words for $w_0$ have diameters
given by the values of $|L_2|$ shown in Table~\ref{L2-size-table}.
\end{theorem}
Unfortunately, types $A_{n-1}, B_n,$ and $I_2(m)$
are the {\it only} irreducible real 
reflection groups\footnote{Note that taking 
{\it products} of hyperplane arrangements
\cite[Def. 2.13]{OrlikTerao} which are supersolvable
preserves supersolvability, and every finite real reflection
group has its reflection arrangement equal to a product
of reflection arrangements for {\it irreducible} real reflection groups.}
whose reflection arrangements are supersolvable;
see Barcelo and Ihrig \cite[Theorem 5.1]{BarceloIhrig}.

\begin{table}[h]
\caption{} \label{L2-size-table}
\newcommand{\str}{\rule[-1.2ex]{0pt}{3.6ex}}
\begin{tabular}{|c|c|c|}\hline
$W$ \str& $|L_2|$ & Does $G_2$ have diameter $|L_2|$? \\
\hline\hline
$A_{n-1}$ \str&  $\frac{1}{24} n(n-1)(n-2)(3n-5)$
   & Yes, by Theorem~\ref{thm:supersolvable-theorem}. \\
\hline
$B_n$     \str&  $\frac{1}{6} n(n-1)(3n^2-5n+1)$
   & Yes, by Theorem~\ref{thm:supersolvable-theorem}. \\
\hline
$D_n$     \str&  $\frac{1}{6} n(n-1)(3n^2 - 11 n + 13)$
   & Yes, for $n \leq 4$; unknown generally. \\
\hline
$E_8$     \str&  4900 & Unknown.\\
\hline
$E_7$     \str&  1281 & Unknown.\\
\hline
$E_6$     \str&  390 & Unknown.\\
\hline
$F_4$     \str&  122 & Unknown.\\
\hline
$H_4$     \str&  722 & Unknown.\\
\hline
$H_3$     \str&  31 & Yes, by Theorem~\ref{thm:Cordovil-result}.\\
\hline
$I_2(m)$     \str&  1& Yes, trivially.\\
\hline
\end{tabular}
\end{table}

We remark here on some of the data related to
Table~\ref{L2-size-table}.  
For type $D_4$, computer calculations\footnote{
These computations give $|\RRR(w_0)|$ in type $D_5$ as $12985968$, and in type $D_6$ as 
$3705762080$.}
by Rob Edman show that 
one has $|\RRR(w_0)|=2316$, in agreement with \cite[\S 7]{Stanley-words}.  In addition,
these calculations exhibit many reduced words for $w_0$ 
giving rise to nodes of the graph $G_2$ which
are {\it $L_2$-accessible} in the sense of 
Definition~\ref{accessibility-definition} 
and Proposition~\ref{prop:accessible-vertex-gives-diameter}
below, showing that the diameter is $|L_2|$.
Interestingly, none of the $L_2$-accessible nodes come from
reduced words that are {\it lexicographically first}
among all reduced words, no matter how one linearly
orders the Coxeter generators.
For $F_4$, these calculations show that
$|\RRR(w_0)|=2144892$, making the full diameter calculation harder,
but again, the computer has checked that none of
the lexicographically first reduced words gives
an $L_2$-accessible node.


The remainder of the paper is structured as follows.
Section~\ref{section:first-two-graphs} establishes formal
definitions for the graph $G_2$ to be studied,
remarking on its connectivity,
its diameter in low dimension, as
well as its relation to mononotone path zonotopes.
Section~\ref{section:distances} 
introduces the notion of a {\it set-valued metric} on a graph,
which is then applied in Section~\ref{section:supersolvable} to prove 
Theorem~\ref{thm:supersolvable-theorem}.
Section~\ref{section:subgraphs} explains how some of these
results adapt for the graphs $G(w)$
when $w$ is not the longest element $w_0$.

\section{Arrangements and the graphs $G_1, G_2$}
\label{section:first-two-graphs}

We review here some of the theory of hyperplane arrangements; see
Orlik and Terao~\cite{OrlikTerao} and Stanley
\cite{Stanley-arrangements} for good references.

As in the Introduction, $\AAA=\{H_1,\ldots,H_N\}$ will
be an arrangement of hyperplanes in $\RR^d$, which is \emph{central} and
\emph{essential}, that is $\bigcap_{i=1}^N H_i = \{0\}$.
The \emph{intersection poset} $L$ for $\AAA$ is the collection
of intersection subspaces $X=\bigcap_{i \in I} H_i$ of subsets
of the hyperplanes, ordered by \emph{reverse} inclusion.  This
makes $L$ a \emph{geometric lattice} (see \cite[\S2.1]{OrlikTerao}),
and therefore {\it graded} or {\it ranked}.  
Let $L= \bigsqcup_{i=0}^d L_i$ be the decomposition into its ranks,
so that 
\begin{itemize}
\item the bottom rank $L_0$ contains only the empty intersection $\RR^d$ itself,
\item the set of atoms $L_1$ is the set of
hyperplanes $\{H_1,\ldots,H_n\}=\AAA$, and
\item the top rank $L_d$ contains only the zero subspace $\{0\}$.
\end{itemize}
The complement $\RR^d \setminus \AAA$ decomposes into
a collection $\CCC$ of connected components called \emph{chambers}.
Given two chambers $c,c'$, define their separation set
$$
L_1(c,c'):=\{ H \in L_1=\AAA: H\text{ separates }c\text{ from }c' \}.
$$

\begin{defn}[The graph $G_1$]
\label{def:first-graph}
Given an arrangement $\AAA$, define a graph $G_1$ whose vertex set is
the set of chambers $\CCC$, and having an edge between two chambers $\{c,c'\}$ exactly
when $|L_1(c,c')|=1$.
\end{defn}

A \emph{(minimal) gallery} from chamber $c$ to chamber $c'$
is a \emph{geodesic} (shortest path)
$$
c:=c_0,c_1,\ldots,c_{d-1},c_d:=c'
$$
in this graph $G_1$.
Fixing one particular choice of a \emph{base chamber} $c_0$,
let $\RRR$ denote the set of all minimal galleries $r$ from
$c_0$ to $-c_0$.

\begin{figure}
\centerline{\includegraphics[width=2in]{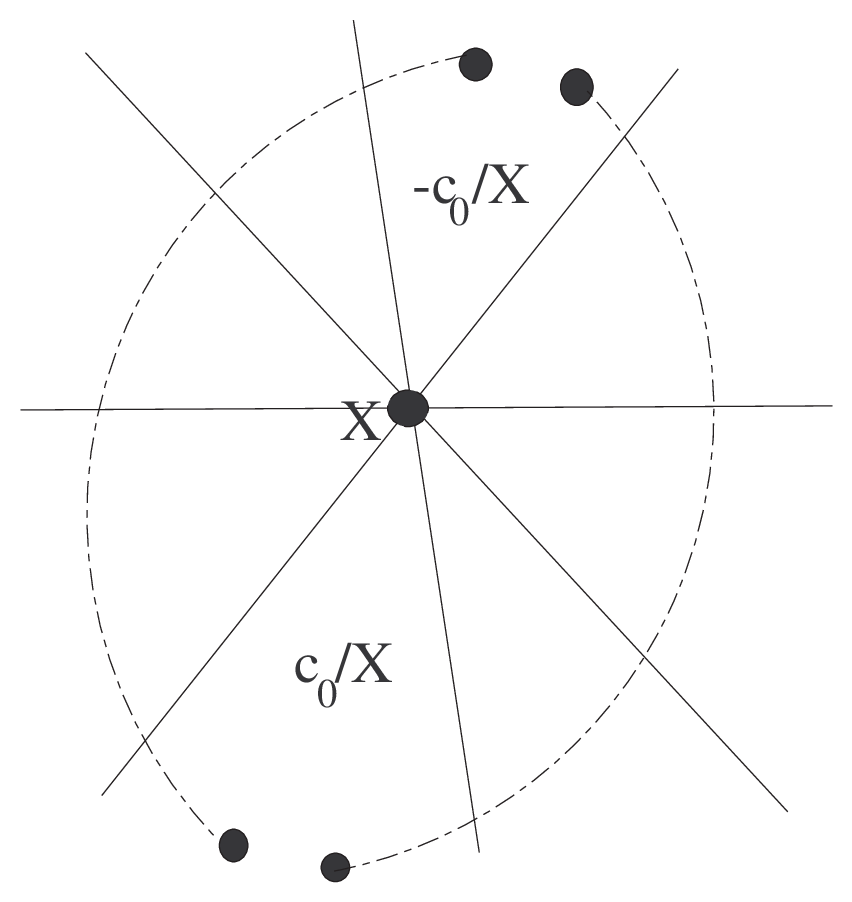}}
\caption{The two possibilities for a minimal gallery $r/X$ from
$c_0/X$ to $-c_0/X$ through $\AAA_X$
when $X$ has codimension-two, each shown as a dashed path.}
\label{fig:rank-two-figure}
\end{figure}

We wish to discuss how the codimension-two intersection subspaces
in $L_2$ can \emph{separate} minimal galleries.
A minimal gallery $r$ in $\RRR$ must cross \emph{every}
hyperplane $H_1,\ldots,H_N$ of $\AAA$ exactly once, and is completely determined by the linear
order in which they are crossed.
Given any intersection subspace $X$, one defines the \emph{localized arrangement}
of hyperplanes in the quotient space $\RR^d/X$
$$
\AAA_X := \{H/X: H \in \AAA \text{ and }H \supseteq X \}.
$$
\noindent
Note that the intersection lattice for $\AAA_X$ may be identified with
the lower interval $[\RR^d,X]$ within $L$.  
For each chamber $c$ of $\AAA$, 
there is a unique chamber $c/X$ of $\AAA_X$ that
contains all additive cosets of the subspace $X$ represented by points of $c$.
A minimal gallery $r$ in $\RRR$ from $c_0$ to $-c_0$ induces a minimal gallery
$r/X$ from $c_0/X$ to $-c_0/X$ in $\AAA_X$.  In particular, when $X$ has
codimension two, so that $\AAA_X$ is an arrangement of lines through the 
origin in the $2$-dimensional plane 
$\RR^d/X$, every minimal gallery $r$ has exactly two possibilities
for the induced minimal gallery $r/X$ from $c_0/X$ to $-c_0/X$ in $\AAA_X$;
see Figure~\ref{fig:rank-two-figure}.

Given two minimal galleries $r,r'$ in $\RRR$
and a codimension-two subspace $X$ in $L_2$, say that
$X$ \emph{separates} $r$ from $r'$ if $r/X \neq r'/X$.  Define
their separation set
$$
L_2(r,r'):=\{ X \in L_2: X\text{ separates }r\text{ from }r' \}.
$$

\begin{defn}[The graph $G_2$]
\label{def:second-graph}
Given the arrangement $\AAA$ and the chosen base chamber $c_0$ in $\CCC$,
define a graph $G_2$ whose vertex set is
the set $\RRR$ of minimal galleries $c_0$ to $-c_0$,
and having an edge between two galleries $\{r,r'\}$ exactly when $|L_2(r,r')|=1$.
\end{defn}

\begin{figure}
\centerline{\begin{overpic}[width=.8\textwidth]{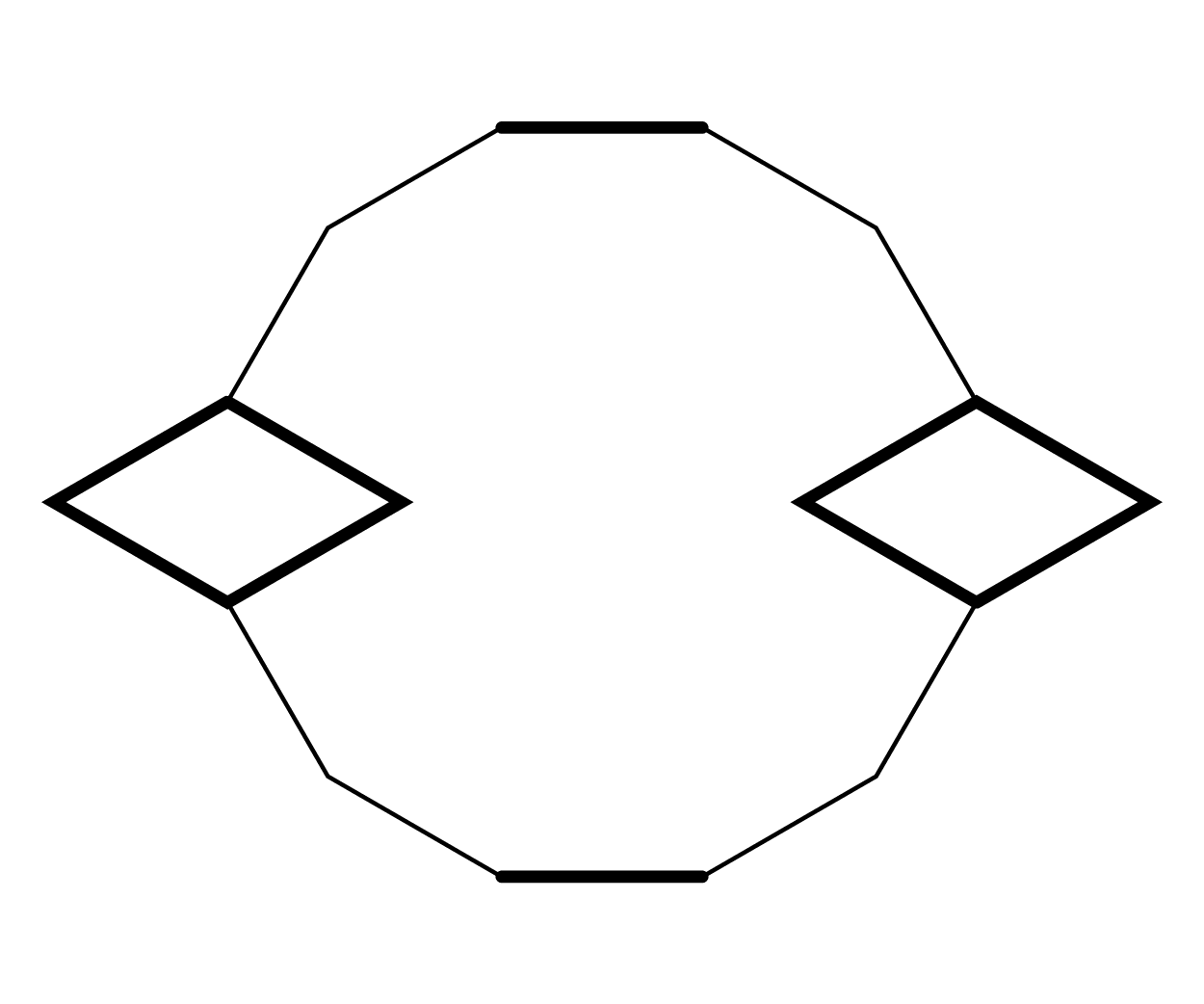}
\put(47,7){$X_{14,23}$} \put(47,75){$X_{14,23}$}
\put(28,70){$X_{234}$} \put(66,12){$X_{234}$}
\put(29,12){$X_{123}$} \put(65,70){$X_{123}$}
\put(16,59){$X_{134}$} \put(77,23){$X_{134}$}
\put(17,23){$X_{124}$} \put(77,59){$X_{124}$}
\put(26,48){$X_{13,24}$} \put(89,35){$X_{13,24}$}
\put(27,35){$X_{12,34}$} \put(88,48){$X_{12,34}$}
\put(4,48){$X_{12,34}$} \put(66,35){$X_{12,34}$}
\put(4,35){$X_{13,24}$} \put(66,48){$X_{13,24}$}
\end{overpic}}
\caption{Labelling the edges $\{r,r'\}$ in $G(w_0)$ for $w_0$
by the unique element $X$ in $L_2(r,r')$}.\label{fig:S4-separation-sets}
\end{figure}

The following proposition points out how separation sets
$L_i(-,-)$ for $i=1,2$ encode the chambers $\CCC$ and galleries $\RRR$.

\begin{proposition}
\label{prop:separation-sets-encode}
Given any fixed base chamber $c_0$, the set $L_1(c_0,c)$ determines the chamber $c$ uniquely.
Given any fixed base gallery $r_0\in \RRR$ from $c_0$ to $-c_0$,
the set $L_2(r_0,r)$ determines the gallery $r$ uniquely.
\end{proposition}
\begin{proof}
The first assertion is clear (and implicit in the discussion of
\cite[\S I]{Edelman}) since $L_1(c_0,c)$
determines on which side of each hyperplane $H$ of $\AAA$
the chamber $c$ lies.

For the second assertion, as noted earlier,
since $L_1(c_0,-c_0) = L_1 = \AAA$, the gallery $r$ from $c_0$ to $-c_0$
must cross \emph{every} hyperplane $H$ of $\AAA$, and $r$ is determined by
the linear order in which it crosses these hyperplanes.  This linear order is
determined by knowing for each pair $H, H'$ which of the two is crossed first.
The latter is determined from the order in which $r$ crosses the hyperplanes of the localized
arrangement $\AAA_X$ for the codimension-two subspace $X:=H \cap H'$, and this is
encoded by the separation set $L_2(r_0,r)$.
\end{proof}

\begin{example}
\label{first-A-example}
The reflection arrangement of type $A_{n-1}$, corresponding
to the symmetric group $W=\sym_n$, has ambient space isomorphic to
$\RR^d$ for $d=n-1$; one identifies $\RR^{n-1}$ 
with the quotient of $\RR^n$ (having coordinates $x_1,\ldots,x_n$)
by the subspace $x_1 = x_2 = \cdots= x_n$.  Its hyperplanes are
$H_{ij}:=\{x_i = x_j\}$ for $1 \leq i < j \leq n$, and its
codimension-two intersection subspaces in $L_2$ are either
of type (i) $X_{ij,k\ell}:=\{x_i=x_j,x_k=x_\ell\}$ or
of type (ii) $X_{ijk}:=\{x_i=x_j=x_k\}$,
corresponding to the braid relations of types (i), (ii) from the Introduction.
The graph $G(w_0)$ from \eqref{fig:S4-example-figure} is redrawn in
Figure~\ref{fig:S4-separation-sets}, with each edge $\{r,r'\}$ labeled
by the unique codimension-two subspace $X$ separating $r$ from $r'$.
\end{example}

We close this section with three remarks on the graph $G_2$.
All of these can be safely skipped by the reader solely interested
in the proof of Theorem~\ref{thm:supersolvable-theorem}.

\begin{remark} (on the connectivity of $G_2$)
\label{rmk:connectivity-credit}
It is not obvious that the graph~$G_2$ is connected
for every real hyperplane arrangement~$\AAA$
and every choice of base chamber $c_0$.
However, as mentioned in the Introduction,
this connectivity of $G_2$ was proven at
the following successively stronger levels of generality:
\begin{itemize}
\item for real reflection arrangements by Tits~\cite{Tits-words},
\item for real simplicial arrangements by Deligne~\cite{Deligne},
\item for all real arrangements by Salvetti~\cite{Salvetti}, and
\item for oriented matroids by Cordovil and Moreira~\cite{CordovilMoreira}.
\end{itemize}

\end{remark}

\begin{remark}(on the diameter of $G_2$ in low dimension)
\label{remark:low-dimension}

When the arrangement $\AAA$ lives in $\RR^d$ for $d \leq 2$,
regardless of the choice of base chamber $c_0$,
the graph $G_2$ is trivial, consisting of a single vertex for $d=1$,
and consisting of two vertices connected by a single edge for $d=2$.

When $d=3$, regardless of the choice of base chamber $c_0$,
the diameter of $G_2$ is exactly $|L_2|$ by the following 
result of Cordovil \cite[Theorem 2.5]{CordovilMoreira} (cited
there as being implicit in \cite[Theorem 2.1]{Cordovil}), 
and proven even more generally for rank $3$ oriented matroids:

\begin{theorem}(Cordovil)
\label{thm:Cordovil-result}
In a $3$-dimensional real central hyperplane arrangement $\AAA$
(or even a rank $3$ oriented matroid), any two minimal 
galleries $r,r'$ from a chamber $c_0$ to its opposite
$-c_0$ can be connected by a 
sequence of at most $|L_2(\AAA)|$ elementary deformations.
\end{theorem}

When $d=4$, this assertion {\it fails} at the level of
generality of oriented matroids;  the authors thank Jim Lawrence
for pointing out how this follows from an important counterexample of 
J. Richter-Gebert, which we recapitulate here; see \cite{OMbook}
for most of the oriented matroid terminology left undefined.  

  A crucial notion is that of a {\it strong map} 
$N \rightarrow M$, where $N,M$ are oriented matroids on the same ground set $E$;
this is defined \cite[\S 7.7]{OMbook}
by requiring that the {\it covectors} of $M$, as a subset of $\{0,+,-\}^E$,
form a subset of the covectors of $N$.  This combinatorially
abstracts the arrangement picture as follows.  When 
$N$ comes from a collection of vectors $\{v_e\}_{e \in E}$
in $\RR^d$, thought of as the normal vectors to the hyperplanes of $\AAA$,
then its covectors are the sign vectors $(\sgn f(v_i))_{i \in E}$ attained
when varying over all linear functionals $f$ in $(\RR^d)^*$;
equivalently, they index the (relatively open) cones of all dimensions
in the decomposition of $\RR^d$ by the hyperplanes of $\AAA$. 
Then a strong map $N \rightarrow M$ abstracts the situation where
$M$ comes from the  image vectors $\{\varphi(v_e)\}_{e \in E}$
under some linear map $\RR^d \overset{\varphi}{\to} \RR^{d'}$.

Richter-Gebert constructs
in \cite[\S3]{RichterGebert} a certain rank $4$ non-realizable, non-Euclidean 
oriented matroid that he calls $R(12)$, having $12$ pseudohyperplanes in general
position.  He shows \cite[Corollary 3.5]{RichterGebert}
that there is a strong map from $R(12) \rightarrow M$ where
$M$ is a uniform rank $2$ oriented matroid, such that the topes
(maximal covectors) of $M$ thought of as a subset of the topes 
of $R(12)$ cannot be contained in the topes of any pseudohyperplane that
extends $R(12)$ by a single element (disproving a conjecture of 
M. Las Vergnas; see \cite[Corollary 3.5]{RichterGebert}).  

Now pick $c_0$ to be any tope of $M$, which is necessarily also
a tope of $R(12)$, and pick $r,-r$ to be the
two unique minimal galleries from $c_0$ to $-c_0$ 
passing through topes of $M$, that is, $r, -r$
pass through the $12$ hyperplanes of $M$ or $R(12)$ in exactly reversed
orders.  Hence if there existed a sequence 
of $\binom{12}{2}=|L_2(R(12))|$ elementary deformations 
connecting $r$ to $-r$, this would lead to a {\it simple allowable
sequence} of permutations of length $12$ in the sense
of Goodman and Pollack;  see \cite[Chapter 6]{OMbook}.  
Such an allowable sequence would then give rise to the topes of
a uniform rank $3$ oriented matroid, containing the topes of $M$, and
coming from a pseudohyperlane that extends $R(12)$ by a single element, 
contradicting Richter-Gebert's result.
\end{remark}

\begin{remark}
(on the relation to monotone path zonotopes)
\label{remark:monotone-path-zonotopes}

We explain here how Billera and Sturmfels' theory of \emph{fiber polytopes}
\cite{BilleraSturmfels} offers an enlightening
perspective on the graph $G_2$, implying good behavior for
certain of its subgraphs.  
%
The reader is referred
to Ziegler \cite[Lectures 7 and 9]{Ziegler-polytopes} for definitions
and terminology omitted in this discussion.

Consider the (central, essential)
arrangement $\AAA$ in $\RR^d$ as the \emph{normal fan} for the \emph{zonotope}
$Z(\AAA)$ which is generated by functionals $\alpha_H\in (\RR^d)^*$
that cut out the hyperplanes $H$, that is, it is the
\emph{Minkowski sum}
$
Z(\AAA) = \sum_{H \in \AAA} [-\alpha_H, \alpha_H]
$
of the line segments $[-\alpha_H, \alpha_H]$.
Then the graph $G_1$ defined in Definition~\ref{def:first-graph}
is exactly the \emph{$1$-skeleton} of this zonotope $Z(\AAA)$.

Now assume that the functionals $\alpha_H$ have been chosen to be
positive on points in the chosen base chamber $c_0$ of $\AAA$.
Then any point $f$ in (the interior of) $-c_0$ gives
a linear functional that achieves its minimum, maximum
values on $Z(\AAA)$ at the vertices whose normal cones are the chambers $c_0, -c_0$.
Let $I=f(Z(\AAA))$ be the interval inside $\RR$ which is the image
of $Z(\AAA)$ under this functional $f$.  Then the \emph{fiber polytope/monotone path polytope}
$$
Z_2 := \Sigma\bigl( Z(\AAA) \overset{f}{\rightarrow} I \bigr) = \Sigma_{f}(Z(\AAA))
$$
discussed by Billera and Sturmfels
in \cite[Theorem 5.3]{BilleraSturmfels} is a $(d-1)$-dimensional
polytope with several interesting properties.

The $1$-skeleton of $Z_2$ turns out to be
a certain subgraph of the graph $G_2$ from
%
%
Definition~\ref{def:second-graph}.
Specifically, minimal galleries $c_0$ to $-c_0$ correspond
to \emph{$f$-monotone paths} $\gamma$ in the $1$-skeleton of
$Z(\AAA)$.  The vertices of $Z_2$ correspond to the subset of
$f$-monotone paths $\gamma$ that are \emph{coherent}, in the sense
that there exists some linear functional $g$ whose
maximum over each fiber $f^{-1}(x) \cap Z(\AAA)$ for $x\in I$
is achieved uniquely at the point $f^{-1}(x)\cap \gamma$;
see \cite[Theorem 2.1]{BilleraSturmfels}.

Furthermore, \cite[Theorem 2.4, Theorem 4.1]{BilleraSturmfels}
imply that $Z_2$ is a \emph{zonotope}, generated by the vectors
$\{ v_{H,H'} \}_{H \neq H' \in \AAA}$
where
$
v_{H,H'}:=f(\alpha_H) \alpha_{H'} - f(\alpha_{H'}) \alpha_H.
$
One can check that, for {\it generic choices} of $f$ within 
the interior of $-c_0$,
any two such generating vectors $v_{H_1,H_2}$ and $v_{H_3,H_4}$
for $Z_2$ are scalar multiples of each other exactly when the
codimension-two intersection subspaces $H_1 \cap H_2$ and $H_3 \cap H_4$
are the same subspace $X$ in $L_2$.

Hence $Z_2$ will then be a zonotope
having exactly $|L_2|$ distinct parallelism classes among its
generating vectors, and the $1$-skeleton 
of $Z_2$ will be a geometrically-distinguished
%
%
subgraph of $G_2$ having the expected diameter $|L_2|$.

%

\end{remark}

\section{Set-valued metrics on graphs}
\label{section:distances}

We introduce some easy observations that apply to the question
of diameter for the graphs $G_1, G_2$ defined in the previous section.

\begin{defn}
Let $G=(V,E)$ be a simple graph on vertex set $V$, meaning that $E$ is a
set of unordered pairs $\{x,y\}$ with $x \neq y\in V$.

The graph-theoretic \emph{distance} $d_G(x,y)$ is the minimum length $d$ of
a path
\begin{equation}
\label{typical-path}
x=v_0,v_1,\ldots,v_{d-1},v_d=y
\end{equation}
with $\{v_i,v_{i+1}\}\in E$ for each~$i$.  Call such a shortest path a {\it geodesic}.

The \emph{diameter} of $G$ is the maximum value of $d_G(x,y)$ over all
$x,y\in V$.
\end{defn}

Note that $d_G(-,-)$ satisfies the usual properties of a \emph{metric} on $V$,
that is,
\begin{itemize}
\item $d_G(x,x)=0$,
\item $d_G(x,y)=d_G(y,x)$, and
\item $d_G(x,z) \leq d_G(x,y)+d_G(y,z)$.
\end{itemize}
We also make the trivial observation that if
$\alpha: V \rightarrow V$ is a \emph{graph automorphism}, meaning a bijection
such that for every edge $\{x,y\}$ in $E$,
the image $\{\alpha(x),\alpha(y)\}$ is also in $E$,
then $\alpha$ takes geodesics to geodesics and preserves distances:
$$
d_G(\,\, \alpha(x) \,\, , \,\, \alpha(y) \,\,) = d_G(x,y).
$$

\begin{defn}
\label{def:set-valued-metric}
For a connected simple graph $G=(V,E)$ and a set $\Omega$, say that a function
\begin{align*}
\Omega(-,-):  V \times V & \longrightarrow  2^\Omega \\
              (x,y)      & \longmapsto      \Omega(x,y)
\end{align*}
is a \emph{set-valued metric  on $G$} if

\begin{itemize}
%
\item[(a)] $\Omega(x,y)=\Omega(y,x)$,
\item[(b)] whenever $\{x,y\}$ is an edge in $E$, one has $|\Omega(x,y)|=1$, and
\item[(c)] $\Omega(x,z) = \Omega(x,y) \symmdiff \Omega(y,z)$
\end{itemize}
where here
$$
A \symmdiff B:=(A \setminus B) \sqcup (B \setminus A)
$$
denotes the symmetric difference of sets.  In particular, the first and third conditions
imply that $\Omega(x,x)=\varnothing$ for any $x$ in $V$.

Here is an equivalent rephrasing of a set-valued metric $\Omega(-,-)$ on $G$:
it is a labelling $\Omega(x,y)$
of each edge $e=\{x,y\}$ in $E$ with an element of $\Omega$ 
in such a way that when one traverses any closed path of edges in the graph, 
each label appears an even number of times.  For {\it any} pair of vertices 
$x,y$ in $V$, not necessarily connected by an edge,
one defines $\Omega(x,y)$ to be the set of labels that appear
an odd number of times on any path from $x$ to $y$.
\end{defn}

\begin{example}
\label{arrangement-example-1}
Given a real hyperplane arrangement $\AAA$, and the graphs $G_i$ for $i=1,2$ defined in
Definitions~\ref{def:first-graph} and \ref{def:second-graph}, one can
easily check that the function $L_i(-,-)$ for $i=1,2$ gives a set-valued
metric.
\end{example}

We begin with two observations about set-valued metrics.

\begin{proposition}
\label{prop:bipartiteness}
A connected simple graph $G=(V,E)$ supports at least one 
set-valued metric if and only if $G$ bipartite.
\end{proposition}
\begin{proof}
Given a set-valued metric $\Omega(-,-)$ on $G$,
choosing any vertex $x_0$ in $V$, one has that $G$ is bipartite
with vertex bipartition $V=V_0 \sqcup V_1$ in which
$$
V_i:= \{y \in V: |\Omega(x_0,y)| \equiv i \mod 2\}. 
$$
Conversely, for $G$ bipartite with vertex bipartition $V=V_0 \sqcup V_1$, 
one can define a trivial set-valued metric 
$\Omega(-,-):  V \times V \longrightarrow  2^{\{e\}}$,
where $\{e\}$ is a singleton, via
$$
\Omega(x,y):=
\begin{cases}
\varnothing &\text{ if }x,y \in V_0\text{ or }x,y \in V_1\\
\{e\} &\text{ otherwise.}
\end{cases}
$$
\end{proof}

\begin{proposition}
\label{prop:set-lower-bound}
A simple graph $G=(V,E)$ with a set-valued metric $\Omega(-,-)$,
has $d_G(x,y) \geq |\Omega(x,y)|$ for all $x,y\in V$.
\end{proposition}
\begin{proof}
A path of length $d$ in~$G$ from $x$ to $y$ as in \eqref{typical-path} leads
to a path
$$
\varnothing = \Omega(x,v_0), \,\,
\Omega(x,v_1), \,\, \ldots, \,\,
\Omega(x,v_{d-1}), \,\,
\Omega(x,v_d) = \Omega(x,y),
$$
where each pair of sets $\Omega(x,v_i), \Omega(x,v_{i+1})$ differs in one element,
namely the unique element of $\Omega(v_i,v_{i+1})$.  Thus $d \geq |\Omega(x,y)|$.
\end{proof}

\begin{defn}
For a set-valued metric $\Omega(-,-)$ on a
simple graph $G=(V,E)$, say an involution $x \mapsto -x$ on the vertex set
$V$ is {\it ($\ZZ_2$-)equivariant} if 
\begin{equation}
\label{eqn:equivariance}
\Omega(x,-y) = \Omega \setminus \Omega(x,y)
\end{equation}
for all $x, y$ in $V$.  This is equivalent,
by property (b) in Definition~\ref{def:set-valued-metric} 
of set-valued metrics, to requiring
only the special case of \eqref{eqn:equivariance} where $y=x$, that
is, equivariance requires only $\Omega(x,-x) = \Omega$ for all $x$ in $V$.
\end{defn}

\begin{example}
\label{arrangement-example-2}
Continuing Example~\ref{arrangement-example-1}, the graphs $G_i$ for $i=1,2$
endowed with the set-valued metrics $L_i(-,-)$ also have equivariant involutions
derived from the linear map $x \mapsto -x$ on $\RR^d$.  For $G_1$, the involution
sends the chamber $c\in \CCC$ to the antipodal chamber $-c$.  For $G_2$,
the involution sends the minimal gallery $r$
$$
c_0,c_1,\ldots,c_{d-1},c_d:=-c_0
$$
to the minimal gallery $-r$ which visits the antipodes of the same chambers in
the reverse order:
$$
c_0=-c_d,-c_{d-1},\ldots,-c_1,-c_0.
$$
Equivalently, $r$ and $-r$ cross the hyperplanes of $\AAA$ in exactly the {\it opposite}
linear orders. 
For example, in the reflection arrangement of type $A_{n-1}$, one can check
that this involution on galleries sends a reduced word
$r=(s_{i_1},s_{i_2},\ldots,s_{i_{\ell-1}},s_{i_\ell})$ for $w_0$
to the reduced word
$
-r=(s_{n-i_\ell},s_{n-i_{\ell-1}},\ldots,s_{n-i_2},s_{n-i_1}).
$
\end{example}

\begin{proposition}
\label{prop:antipode-lower-bound}
A simple graph $G=(V,E)$ with a set-valued metric $\Omega(-,-)$ and
an equivariant involution always has diameter at least $|\Omega|$.
\end{proposition}
\begin{proof}
By Proposition~\ref{prop:set-lower-bound}
\begin{equation*}
d_G(x,-x) \geq |\Omega(x,-x)| =|\Omega \setminus \Omega(x,x)| = |\Omega|.
\qedhere
\end{equation*}
\end{proof}

\begin{defn}
\label{accessibility-definition}
For a set-valued metric $\Omega(-,-)$ on a
simple graph $G=(V,E)$, say that a vertex $x_0$ in $V$ is \emph{$\Omega$-accessible}
if $d_G(x_0,y) = |\Omega(x_0,y)|$ for every $y\in V$.
\end{defn}

\begin{example}
\label{arrangement-example-3}
Continuing Examples~\ref{arrangement-example-1} and
\ref{arrangement-example-2},  it was observed by
Edelman (see \cite[Proposition 1.1]{Edelman})
that for every real hyperplane arrangement $\AAA$, every vertex
in the graph $G_1$ is $L_1$-accessible:  given any two chambers $c,c'$,
a straight-line path between generic points in $c,c'$ gives
a path of length $|L_1(c,c')|$ between their corresponding vertices of $G_1$.
\end{example}

\begin{example}
\label{inaccessibility-example}
Of the $16$ reduced words for $w_0$ in $W=\sym_4$ shown in
Figure \ref{fig:S4-example-figure}, there are exactly four which \emph{do not}
give $L_2$-accessible vertices for the graph $G_2$ shown, namely
the four words
\begin{equation}
\label{four-inaccessible-words}
\{213213, \quad 231231, \quad 132132, \quad 312312\}.
\end{equation}
We check that none of the four words $r_0$ in this set \eqref{four-inaccessible-words}
is $L_2$-accessible.  Scrutiny of
Figures~\ref{fig:S4-example-figure}
and~\ref{fig:S4-separation-sets}
shows that there are exactly two words
at the maximum distance $7(=|L_2|)$ from such an $r_0$, namely its antipodal word
$-r_0$, and a second word $r \neq -r_0$ having
$d_{G_2}(r_0,r) = 7 > 5 =|L_2(r_0,r)|$.  In particular, any such
pair $\{r_0,r\}$ provides an example that answers negatively
a question of Autord and Dehornoy \cite[Question 1.9]{AutordDehornoy}.

On the other hand, one can check using Figures~\ref{fig:S4-example-figure} 
and~\ref{fig:S4-separation-sets}, via brute force (mitigated by some 
symmetry), that all 12 of the other words $r_0$ are $L_2$-accessible.

\end{example}

Our goal in Section~\ref{section:supersolvable}
will be to show that for supersolvable arrangements
$\AAA$, when one chooses the base chamber $c_0\in \CCC$ incident to a chain of modular
flats, there is a choice of base gallery $r_0\in \RRR$ which is $L_2$-accessible.
Therefore, in this case, the diameter for $G_2$ will be
determined by the next proposition.

\begin{proposition}
\label{prop:accessible-vertex-gives-diameter}
Assume one has a simple graph $G=(V,E)$ with a set-valued metric $\Omega(-,-)$,
and an involution $v \mapsto -v$ on $V$ which is both
equivariant and a graph automorphism of $G$.

If $V$ contains an $\Omega$-accessible vertex $x_0$,
then the diameter of $G$ is exactly $|\Omega|$.
\end{proposition}
\begin{proof}
By Proposition~\ref{prop:antipode-lower-bound}, it suffices to
show that $d_G(x,y) \leq |\Omega|$ for all $x,y$.  This follows from
these equalities and inequalities, justified below:
\begin{align*}
2 d_G(x,y)
&= d_G(x,y) + d_G(x,y) \\
&\overset{(1)}{\leq} d_G(x,x_0) + d_G(x_0,y) + d_G(x,-x_0) + d_G(-x_0,y)  \\
&\overset{(2)}{=} d_G(x,x_0) + d_G(x_0,y) + d_G(-x,x_0)+ d_G(x_0,-y)  \\
&\overset{}{=} \left( d_G(x,x_0) + d_G(-x,x_0) \right) +
 \left(  d_G(x_0,y) + d_G(x_0,-y) \right) \\
&\overset{(3)}{=} \left( |\Omega(x,x_0)| + |\Omega(-x,x_0)| \right) + 
 \left( |\Omega(x_0,y)| + |\Omega(x_0,-y)| \right) \\
&\overset{(4)}{=} |\Omega| + |\Omega| 
= 2|\Omega|
\end{align*}
Inequality (1) twice uses the triangle inequality for the metric $d_G(-,-)$.
Equality (2) twice uses the assumption that $v \mapsto -v$ is a
a graph automorphism.
Equality (3) four times uses the $\Omega$-accessibility of $x_0$.
Equality (4) twice uses the assumption of equivariance.
\end{proof}

Applying Propositions~\ref{prop:antipode-lower-bound} and
\ref{prop:accessible-vertex-gives-diameter} to the graphs $G_1, G_2$
immediately gives the following.

\begin{corollary}
\label{cor:two-graphs-diameters}
For any real hyperplane arrangement $\AAA$, the graph $G_1$ has
diameter {\bf exactly} $|L_1|$.
For any choice of a base chamber $c_0$ of $\AAA$, the graph $G_2$ has
diameter {\bf at least} $|L_2|$.
\end{corollary}

\section{Supersolvable arrangements}
\label{section:supersolvable}

We wish to first review the definition and some
properties of supersolvable arrangements~\cite{Stanley-supersolvable},
\cite[\S 2.3]{OrlikTerao}, \cite[\S 4]{BjornerEdelmanZiegler}, and
then apply this to prove Theorem~\ref{thm:supersolvable-theorem}.

\begin{defn}
\label{supersolvable-defn}
Given a finite geometric lattice $L$ with 
rank function $\rho$, say that $x$ in $L$ a {\it modular element} of $L$ if
$
\rho(x \vee y) +\rho(x \wedge y) = \rho(x)+\rho(y)
$
for all $y$ in $L$.
The lattice $L$ is called \emph{supersolvable} if it contains
an {\it $M$-chain}, that is, a maximal chain of modular elements.
\end{defn}

A hyperplane arrangement $\AAA$ is called {\it supersolvable}
when its intersection lattice $L$ is supersolvable.
An inductive rephrasing, due to Bj\"orner, Edelman and
Ziegler \cite[Thm. 4.3]{BjornerEdelmanZiegler},
will be more useful for our purposes.
Given a real (central, essential) hyperplane arrangement $\AAA$ in $\RR^d$,
an element $\ell$ of $L_{d-1}$ is called a \emph{coatom}.  Thus coatoms $\ell$
are {\it lines} obtained by intersecting the hyperplanes. 

\begin{proposition}\cite[Thm. 4.3]{BjornerEdelmanZiegler}
\label{supersolvability-rephrasing-proposition}
Let $\AAA$ be a hyperplane arrangement.
\begin{enumerate}
\item[(i)]
A coatom $\ell$ of $L$ is a {\it modular} coatom if and only if for 
every pair $H, H'$ of distinct hyperplanes
of $\AAA$ \emph{not} containing $\ell$, there exists a hyperplane $H''$ of
$\AAA$ containing both $\ell$ and $H \cap H'$ (that is to say,
the hyperplane $H''=\ell+H \cap H'$ is in $\AAA$).
\item[(ii)]
A hyperplane arrangement $\AAA$ is supersolvable
if and only it satisfies the following inductive definition: either $\AAA$ has rank
$d=1$, or its intersection lattice $L$ contains a modular coatom $\ell$
for which the localized arrangement $\AAA_\ell$ of rank $d-1$ is supersolvable.
\end{enumerate}
\end{proposition}

\begin{example}
When $d=2$ the arrangement $\AAA$ is always supersolvable, as
any of its hyperplanes (=lines) is a modular coatom.
\end{example}

\begin{example}
\label{type-A-example-1}
Recall from Example~\ref{first-A-example} that the
reflection arrangement $\AAA$ of type $A_{n-1}$, corresponding
to $W=\sym_n$, lives in $\RR^n/\{x_1 = x_2 = \cdots= x_n\}$,
and has hyperplanes $H_{ij}:=\{x_i = x_j\}$.
The line $\ell$ defined by
$\{x_1=x_2= \cdots = x_{n-1}\}$ is a modular coatom for $A_{n-1}$:
any two typical hyperplanes $H_{in}, H_{jn}$ for $i,j < n$ that do not contain $\ell$
will have $\ell + \left( H_{in} \cap H_{jn} \right) = H_{ij}$, which is another
hyperplane in the arrangement $A_{n-1}$.
The localization $\AAA_\ell$ is isomorphic to the reflection arrangement of
type $A_{n-2}$, and hence
one can iterate this construction to show that the arrangement of type $A_{n-1}$ 
is supersolvable.
\end{example}

\begin{example}
\label{type-B-example-1}
The reflection arrangement $\AAA$ of type $B_n$ lives in $\RR^n$,
and consists of all hyperplanes of the form
\begin{align*}
H^+_{ij}&:= \{x_i = + x_j\}, \\
H^-_{ij}&:= \{x_i = - x_j\}, \\
H^-_{ii}&:= \{x_i = 0 \}.
\end{align*}
Let us check that the line $\ell$ 
defined by $\{x_1=x_2= \cdots = x_{n-1}=0\}$
is a modular coatom for $B_n$.  Pairs of
hyperplanes not containing $\ell$ either come from choosing
\begin{enumerate}
\item[$\bullet$]
two indices $i,j \leq n-1$ and 
two signs $\alpha,\beta=\pm 1$, giving
hyperplanes $H^{\alpha}_{in}, H^{\beta}_{jn}$ 
which satisfy
$\ell+ \left( H^{\alpha}_{in} \cap H^{\beta}_{jn} \right)
= H^{\alpha \cdot \beta}_{ij}$, a hyperplane of $B_n$, or 
\item[$\bullet$]
an index $i \leq n-1$ and a sign $\alpha$,
giving hyperplanes $H^{-}_{nn}, H^{\alpha}_{in}$,
which satisfy
$\ell+ \left( H^{-}_{nn} \cap H^{\alpha}_{in} \right)
= H^{-}_{ii}$, a hyperplane of $B_n$.
\end{enumerate}
The localization $\AAA_\ell$ is isomorphic to the reflection arrangement 
of type $B_{n-1}$, and hence
one can iterate this construction to show that the arrangement of type 
$B_n$ is supersolvable.
\end{example}


Given an intersection subspace $X$, define the map $\pi_X$ sending chambers of $\AAA$
to their corresponding chamber in the localization $\AAA_X$:
\begin{align*}
\pi_X: \CCC=:\CCC(\AAA) &\longrightarrow \CCC(\AAA_X) \\
                       c&\longmapsto     c/X
\end{align*}
Say that a chamber $c$ is \emph{incident} to a subspace $X$ if
the closure of $c$ intersects $X$ in a subcone of the same dimension as $X$.

\begin{proposition}
\label{prop:modular-coatom-prop}
Assume that $\ell$ is a modular coatom for $\AAA$, and
$c\in \CCC$ is a chamber incident to $\ell$.
\begin{enumerate}[(i)]
\item (compare with the discussion before Theorem 4.4 of~\cite{BjornerEdelmanZiegler})
There is a linear order on the fiber
$$
\pi^{-1}_\ell(\pi_\ell(c))=\{c_1(=c),c_2,c_3,\ldots,c_t\}
$$
such that the sets $L_1(c,c_i)$ for $i=1,2,\ldots,t$ are nested:
$$
\varnothing = L_1(c,c_1) \subset L_1(c,c_2)
\subset \cdots \subset
L_1(c,c_t) = \AAA \setminus \AAA_\ell,
$$
This induces a linear order $H_1,H_2, \ldots $ on $\AAA \setminus \AAA_\ell$
such that $H_i$ is the unique hyperplane in $L_1(c,c_i) \setminus L_1(c,c_{i-1})$.

\item
Using the linear order 
$H_1,H_2, \ldots $ on $\AAA \setminus \AAA_\ell$
from part (i), if  $i < j < k$ and
if the chamber $c$ incident to $\ell$ is also incident to $\ell+H_i \cap H_k$, then
\begin{equation}
\label{three-intersections-equal}
H_i \cap H_j = H_i \cap H_k = H_j \cap H_k.
\end{equation}
\end{enumerate}
\end{proposition}
\begin{proof}

\noindent
{\sf Proof of assertion (i):}
Assume for the sake of contradiction that
there exist two chambers $c_i, c_j$ with $\pi_\ell(c_i)=\pi_\ell(c_j)=\pi_\ell(c)$
and two hyperplanes $H_i, H_j$ not containing $\ell$ for which
\begin{equation}
\label{incomparable-hyperplanes}
\begin{split}
H_i &\in L_1(c, c_i) \setminus L_1(c, c_j) \\
H_j &\in L_1(c, c_j) \setminus L_1(c, c_i).
\end{split}
\end{equation}
By the modularity of the coatom $\ell$,
the hyperplane $H:=\ell+H_i \cap H_j$ is in $\AAA$.
Consider the codimension-two subspace
$$
X:=H_i \cap H_j = H \cap H_i = H \cap H_j
$$
and the local picture for the lines and chambers
$$
H/X, \; H_i/X, \; H_j/X, \; c/X, \; c_i/X, \; c_j/X
$$
within the localized rank two arrangement $\AAA_X$.
Then \eqref{incomparable-hyperplanes} together with
the assumption that $H$ contains the line $\ell$ incident to $c$
forces this local picture to be as in Figure~\ref{fig:incomparability-figure}(i).
In particular, it forces $H$ to separate $c_i, c_j$,
and since $\ell \subset H$, this contradicts $\pi_\ell(c_i)=\pi_\ell(c_j)$.

\begin{figure}
\centerline{\includegraphics[width=3in]{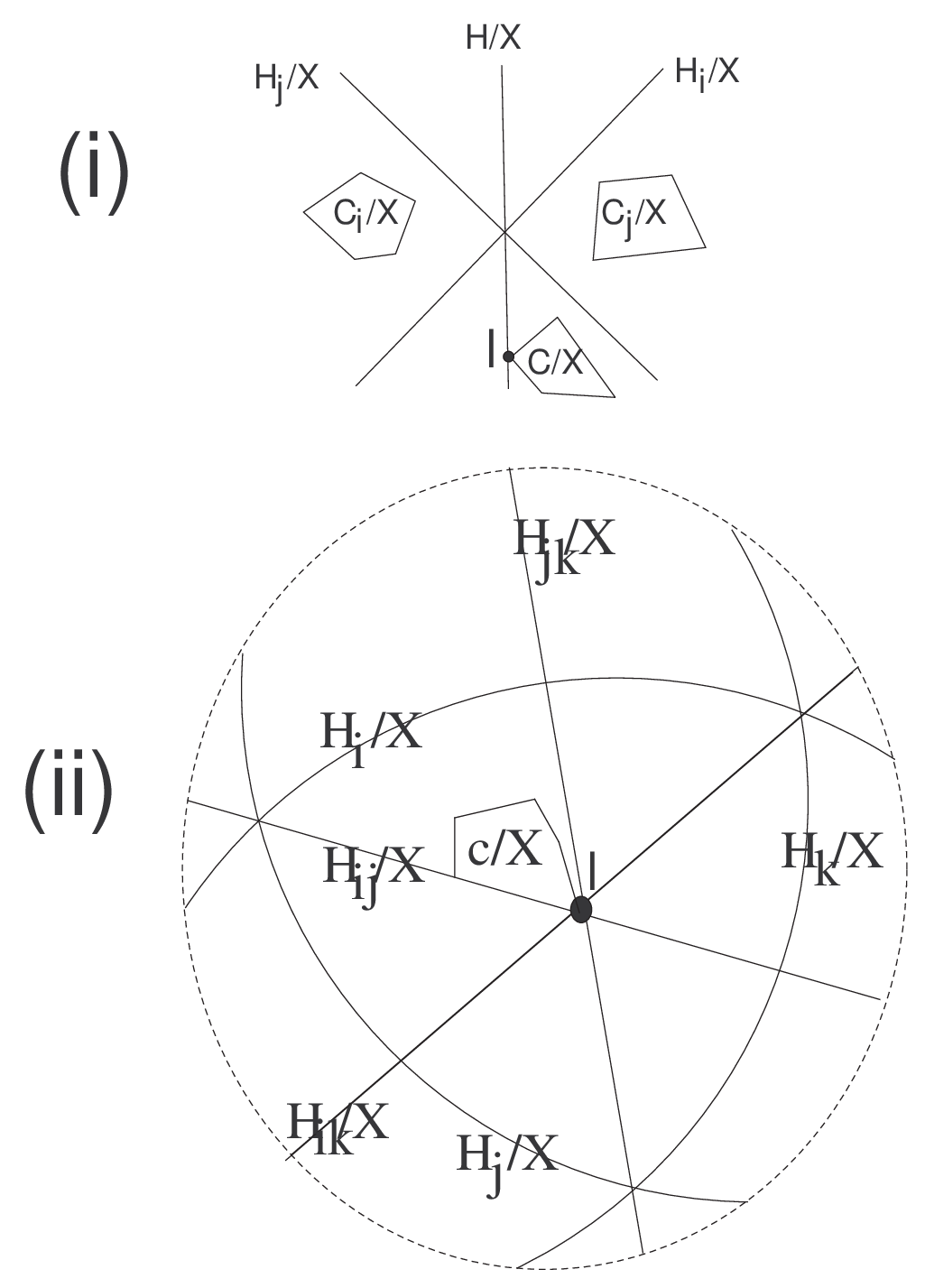}}
\caption{(i) Local picture illustrating why the 
fibers $\pi_\ell^{-1}(\pi_\ell(c))$ are linearly ordered
by inclusion of $L(c,-)$. \newline
(ii) Local picture illustrating why incidence of $c$ to $\ell + H_i \cap H_k$ and $i<j<k$ 
forces the equality \eqref{three-intersections-equal}.
}
\label{fig:incomparability-figure}
\end{figure}

\vskip.1in
\noindent
{\sf Proof of assertion (ii):}
Assume for the sake of contradiction
that $i < j < k$ and $c$ is incident to $\ell+H_i \cap H_k$,
but \eqref{three-intersections-equal} fails,
so that the intersection
$X:=H_i \cap  H_j \cap  H_k$ is of codimension three,
not two.  Note that $\ell$ is contained in none of $H_i, H_j, H_k$.
Therefore since $\ell$ is a modular coatom, each of
the following three hyperplanes containing both $\ell$ and $X$ must also
be a hyperplane in $\AAA$:
\begin{align*}
H_{ij}&:= \ell+H_i \cap H_j \\
H_{ik}&:= \ell+H_i \cap H_k \\
H_{jk}&:= \ell+H_j \cap H_k.
\end{align*}
Now consider the local picture for
$$
H_i/X, \,\,
H_j/X, \,\,
H_k/X, \,\,
H_{ij}/X, \,\,
H_{ik}/X, \,\,
H_{jk}/X, \,\,
\ell/X
$$
within $\AAA_X$, which after an
invertible linear transformation of $\RR^d/X$,
can be made to look as in Figure~\ref{fig:incomparability-figure}(ii).

Recalling that $c$ is incident to the line $\ell$,
the condition $i < j < k$,
forces $c/X$ to be in the chamber shown, so that as one starts in
$c$ and moves away from $\ell$ staying within the same chamber
of $\AAA_\ell$, one crosses the hyperplanes $H_i, H_j, H_k$ in this order;
if $c/X$ lies in any of the other five chambers incident to $\ell/X$
in this figure, one will cross $H_i, H_j, H_k$ in a different order.

However, Figure~\ref{fig:incomparability-figure}(ii) also shows that
this location for $c/X$ contradicts the incidence of
$c$ to $H_{ik}=\ell+H_i \cap H_k$.
\end{proof}

\begin{defn}
A \emph{(maximal) flag} $F=\{X_i\}_{i=0}^d$ is a chain of intersection subspaces in $L$
$$
\{0\}=X_0 \subset X_1 \subset X_2 \subset \cdots \subset X_{d-1} \subset X_d=\RR^d
$$
in which $X_i$ is of dimension $i$.
Say that a chamber $c$ is \emph{incident} to the flag $F$ if $c$ is incident to each of
the $X_i$. If $c_0$ is a chamber incident to a
flag $F$, say that a minimal gallery $r$ from
$c_0$ to $-c_0$ is \emph{incident} to the flag $F$ if it first crosses the
unique hyperplane in $\AAA_{X_{d-1}}$ (namely $X_{d-1}$ itself), then crosses
the hyperplanes in $\AAA_{X_{d-2}} \setminus \AAA_{X_{d-1}}$, etc.,
always crossing the hyperplanes in $\AAA_{X_j}$ before those in
$\AAA_{X_i} \setminus \AAA_{X_j}$ whenever $i<j$.

According to defintion of supersolvability 
(Definition~\ref{supersolvable-defn}) and its rephrasing in
Proposition \ref{supersolvability-rephrasing-proposition},
$\AAA$ is supersolvable if and only if it has an $M$-chain 
or \emph{modular flag}, in which $X_1$ is a modular coatom for $\AAA$,
while $X_2/X_1$ is a modular coatom in the localized
arrangement $\AAA_{X_1}$, and generally $X_i/X_{i-1}$ is a modular
coatom in the localized arrangement $\AAA_{X_{i-1}}$.
\end{defn}

\begin{proposition}
\label{lifting-proposition}
Let $\ell$ be a coatom of the intersection lattice $L$ for
an arrangement $\AAA$. The map
\begin{align*}
\CCC(\AAA)&\overset{\pi_\ell}{\longrightarrow} \CCC(\AAA_\ell) \\
          &c \longmapsto   c/\ell
\end{align*}
becomes $2$-to-$1$ when restricted to the subset of
chambers of $\AAA$ incident to $\ell$.

Consequently\footnote{Compare this with the discussion before Theorem 4.4 of~\cite{BjornerEdelmanZiegler}}, 
when $\AAA$ is a real central essential arrangement in $\RR^d$,
there are exactly $2^d$ chambers incident to each maximal flag $F$.
\end{proposition}
\begin{proof}
The two chambers of a fiber $\pi^{-1}_\ell(\pi_\ell(c))$ that contain $\ell$ are
the two chambers $c^+, c^-$ in the fiber whose closures contain
$\ell^+, \ell^-$, the two rays (half-lines) comprising the line $\ell$.
\end{proof}

\begin{theorem}
\label{thm:supersolvable-accessibility}
Let $\AAA$ be a real (central, essential) hyperplane arrangement in $\RR^d$ which is
supersolvable, let $F:=\{X_i\}_{i=0}^d$ be a modular flag for $\AAA$,
and let $c_0$ be any chamber incident to $F$.

\begin{enumerate}[(i)]
\item
There is a unique minimal gallery
$r_0$ from $c_0$ to $-c_0$ incident to $F$.
\item
This minimal gallery $r_0$ is an $L_2$-accessible vertex
for the graph $G_2$ on the galleries $\RRR$ from $c_0$ to $-c_0$.
\end{enumerate}
\end{theorem}
\begin{proof}

\vskip.1in
\noindent
{\sf Proof of assertion (i).}
Proceed by induction on $d$, with the base case $d=1$ being trivial.
In the inductive step, let $\ell^+$ be the half-line of $\ell$ contained
in the closure of $c_0$.  Then the unique minimal gallery $r_0$ from
$c_0$ to $-c_0$ incident to $F$ is constructed as follows,
in order to have it cross all the hyperplanes in $\AAA_\ell$ first:
\begin{enumerate}[(a)]
\item
Apply induction to the $(d-1)$-dimensional supersolvable arrrangement $\AAA_\ell$,
to find the unique gallery from $c_0/\ell$ to $-c_0/\ell$ incident to $F/\ell$.
\item
Begin the gallery $r_0$ by lifting each chamber
$c/\ell$ of $\AAA_\ell$ in this gallery from (a)
to the unique chamber $c$ in $\AAA$ whose closure contains $\ell^+$.
\item
After going through the chambers in this lifted gallery from (b),
ending in a gallery called $c$ incident to $\ell$,
one now has no choice about how to complete the rest of $r_0$: one
must cross the hyperplanes $\AAA \setminus \AAA_\ell$ in the
linear order given by Proposition~\ref{prop:modular-coatom-prop}(i).
\end{enumerate}

\begin{figure}
\centerline{\includegraphics[width=3in]{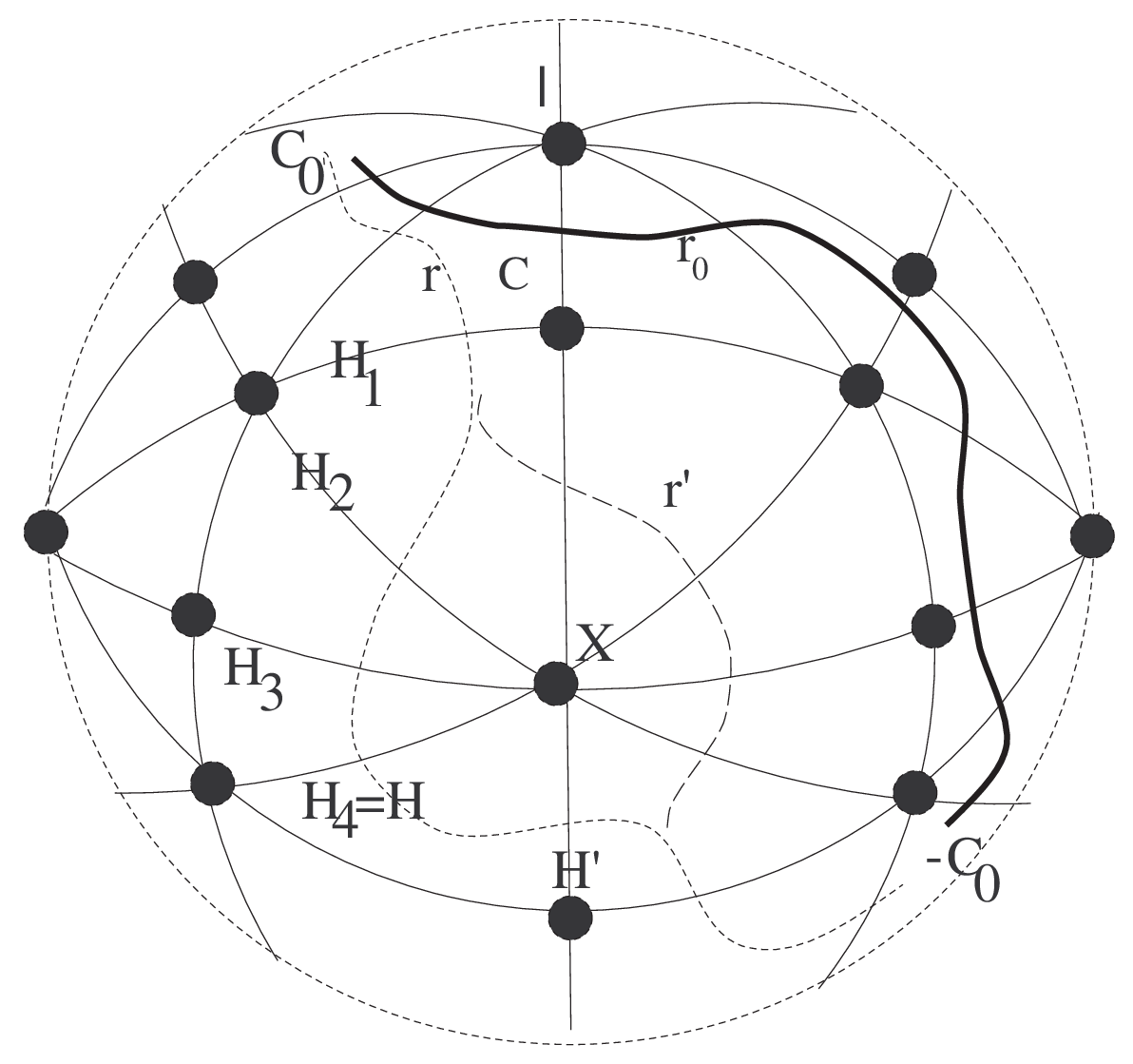}}
\caption{The reflection arrangement of type $B_3$, drawn as
great circles on a unit sphere.  The chamber $c_0$ and
its antipode $-c_0$ are labeled, as well as the three galleries
$r_0, r, r'$ from $c_0$ to $-c_0$ which appear in the proof of
Theorem~\ref{thm:supersolvable-accessibility}}
\label{fig:supersolvable-galleries-figure}
\end{figure}

For example, for the reflection arrangement of
type $A_3$, the base word/gallery $r_0=121321$, indexing the vertex at the bottom left
in Figure~\ref{fig:S4-example-figure}, is the unique gallery $r_0$ incident to the
modular flag described in Example~\ref{type-A-example-1};
this gallery is discussed further below in Example~\ref{type-A-example-2}.
Figure~\ref{fig:supersolvable-galleries-figure} shows
the reflection arrangement of type $B_3$, with
the unique gallery $r_0$ incident to a certain modular flag labeled.

\vskip.2in
\noindent
{\sf Proof of assertion (ii).}
Proceed by induction on $d$, with the base case $d=2$ being trivial.
In the inductive step, it suffices to show that
for any gallery $r \neq r_0$ from $-c_0$ to $c_0$, there
exists another gallery $r'$ having $L_2(r_0,r') \subset L_2(r_0,r)$ and
$|L_2(r,r')|=1$.
There are two cases.

\vskip.1in
\noindent
{\sf Case 1.}
The gallery $r$ crosses all the hyperplanes in $\AAA_\ell$ before
crossing any hyperplanes of $\AAA \setminus \AAA_\ell$.

Then just as with $r_0$, the gallery $r$ 
must cross the hyperplanes in
$\AAA \setminus \AAA_\ell$ in the linear order given
by  Proposition~\ref{prop:modular-coatom-prop}.
Hence the galleries $r$ and $r_0$ can differ only in
their initial segments $\hat{r}, \hat{r_0}$ where they
cross all the hyperplanes of $\AAA_\ell$.  Applying
induction on dimension to the quotient galleries 
$\hat{r}/\ell, \hat{r_0}/\ell$, there exists a gallery $\hat{r'}/\ell$
from $c_0/\ell$ to $-c_0/\ell$ 
in $\AAA_\ell$ having $L_2(\hat{r_0}/\ell,\hat{r'}/\ell) 
\subset L_2(\hat{r_0}/\ell,\hat{r}/\ell)$ and
$|L_2(\hat{r}/\ell,\hat{r'}/\ell)|=1$.
The desired gallery $r'$ from $c_0$ to $-c_0$
is then obtained by first lifting $\hat{r'}/\ell$ as in
Proposition~\ref{lifting-proposition}, and
then completing it by crossing the hyperplanes in
$\AAA \setminus \AAA_\ell$ in the linear order given
by  Proposition~\ref{prop:modular-coatom-prop}.  

\vskip.1in
\noindent
{\sf Case 2.}
The gallery $r$ crosses some hyperplane of $\AAA \setminus \AAA_\ell$
before it has finished crossing all the hyperplanes in $\AAA_\ell$.

Then there must exist at least one ordered pair of hyperplanes
$(H,H')$ crossed consecutively by $r$ that has $H\in \AAA \setminus \AAA_\ell$
and $H'\in \AAA_\ell$.  Find the earliest occurrence\footnote{Compare
the rest of this proof
with the proof of \cite[Lemma 1.2]{AutordDehornoy}, which it generalizes.}
of such a pair $(H,H')$.

Thus $r$ begins by crossing some (possibly empty) sequence of hyperplanes
in $\AAA_\ell$, reaching some chamber $c$ incident to $\ell$,
and then immediately thereafter crosses a sequence of hyperplanes
$H_1,H_2,\ldots,H_{t-1},H_t=:H$
that are all in $\AAA \setminus \AAA_\ell$, before crossing $H'$.
See Figure~\ref{fig:supersolvable-galleries-figure} for
an illustration of a typical gallery $r$, with the chamber $c$ and
the hyperplanes $H_1,H_2,\ldots,H_{t-1},H_t(=:H),H'$
labeled as in this proof;  in this example, $t=4$.

Let $X:=H \cap H'$, a codimension-two subspace in $L_2$.
Note that as $H'=\ell + X$, this $H'$ is the
unique hyperplane in the rank two subarrangement
$\AAA_X$ that also lies in $\AAA_\ell$.
We wish to determine \emph{exactly} when the other hyperplanes
$\AAA_X \setminus \{H'\}$ are crossed by the gallery $r$.

Note that since the quotient gallery $r/\ell$ would cross
$H'/\ell$ to leave $c/\ell$, this chamber $c/\ell$ must be
incident to $H'/\ell$.  As $c$ is
incident to $\ell$, this means that $c$ is incident to $H'$.

Now since $r$ visits the chamber $c$ which is incident to $H'$, but then
crosses another hyperplane $H$ that contains $X$ before crossing $H'$,
it must be that $r$ crosses \emph{every other} hyperplane of the rank two arrangement
$\AAA_X$ before crossing $H'$.

On the other hand, since $r$ only crossed hyperplanes in $\AAA_\ell$ before
reaching $c$, it must be that
$$
\AAA_X \setminus \{H'\} \subset \{H_1,H_2,\ldots,H_{t-1},H_t(=H)\}.
$$

Also note that the hypotheses of Proposition~\ref{prop:modular-coatom-prop}(ii)
are satisfied by $c, \ell$ and by any two hyperplanes $H_i, H_k$
lying in $\AAA_X \setminus \{H'\}$, since $c$ is incident to
$$
H'=\ell+ X = \ell+H_i \cap H_k.
$$
This means that for any $j$ with
$i < j < k$ one has $H_j$ also in  $\AAA_X \setminus \{H'\}$.
Thus the hyperplanes of $\AAA_X \setminus \{H'\}$ occur
\emph{consecutively} within the list $(H_1,\ldots,H_{t-1},H_t)$.
That is, we have shown that there is some index $s \leq t$ for which
$$
\{ H_s,H_{s+1},\ldots,H_{t-1},H_t(=H),H' \} = \AAA_X.
$$

Now let $r'$ be the gallery obtained from $r$ by reversing this
consecutive sequence of crossings\footnote{The fact that another such gallery $r'$ exists turns out to be true, but was
  not justified here nor in the journal version of this paper.
  The authors thank Thomas McConville for pointing out this gap, and for fixing it-- see Section 6 of his preprint
``Biclosed sets in real hyperplane arrangements'' {\tt arXiv:1411.1305}.} $(H_s,H_{s+1},\ldots,H_{t-1},H_t(=H),H')$
of the hyperplanes in $\AAA_X$.
By construction, $L_2(r,r') = \{ X \}$.

We must also check that $L_2(r_0,r') \subset L_2(r_0,r)$,
i.e. that $X$ is not in $L_2(r_0,r')$.
This follows because $r_0$ is incident to $F$, so it must
cross the hyperplane $H'\in \AAA_\ell$ before it can cross
the hyperplane $H\in \AAA \setminus \AAA_\ell$.  Thus in
regard to its order of crossing the hyperplanes of
$\AAA_X$, the gallery $r_0$ agrees with $r'$, not with $r$.
\end{proof}

Theorem~\ref{thm:supersolvable-theorem} is now immediate from
Theorem~\ref{thm:supersolvable-accessibility}(ii) and
Proposition~\ref{prop:accessible-vertex-gives-diameter}.

\begin{example}
\label{type-A-example-2}
Continuing
Example~\ref{type-A-example-1}, for the arrangement $\AAA$ of
type $A_{n-1}$, one can choose as modular flag
$F:=\{X_i\}_{i=0}^{n-1}$ where
$$
X_i:=\{x_1=x_2=\cdots=x_{n-i} \}.
$$
The chambers $c_w\in \CCC$ may be indexed by
permutations $w$
in $\sym_n$, with defining inequalities
$x_{w(1)} < \cdots < x_{w(n)}$.
The chamber $c_0$ corresponding to the identity permutation is incident to
the above modular flag $F$.  The unique gallery $r_0$ from $c_0$ to
$-c_0$ incident to $F$ crosses the hyperplanes in this order:
\begin{align*}
&H_{12},\\
&H_{13}, H_{23}, \\
&H_{14}, H_{24}, H_{34}, \\
&H_{15}, H_{25}, H_{35}, H_{45}, \ldots
\end{align*}
Using the Coxeter generators $S=\{s_1,\ldots,s_{n-1}\}$
for $W=\sym_n$, in which $s_i$ is the adjacent transposition $(i,i+1)$,
this gallery $r_0$ corresponds to the following
reduced decomposition for $w_0$:
\begin{align*}
&(s_1,  \\
&s_2, s_1, \\
&s_3, s_2, s_1, \\
&s_4, s_3, s_2, s_1, \ldots ).
\end{align*}
\end{example}

\begin{example}
\label{type-B-example-2}
Continuing Example~\ref{type-B-example-1}, 
for the arrangement $\AAA$ of type $B_n$, one can choose
as modular flag $F:=\{X_i\}_{i=0}^n$ where
$$
X_i:=
\{x_1=x_2=\cdots=x_{n-i}=0\}.
$$
The chambers $c_w\in \CCC$ may be indexed by
signed permutations $w$ with defining inequalities
$0< \epsilon_1 x_{w(1)} < \cdots < \epsilon_n x_{w(n)}$
if $w$ sends the standard basis vector $e_j$ to $\epsilon_j e_{w(j)}$
with $\epsilon_j \in \{\pm 1\}$.
The chamber $c_0$
corresponding to the identity permutation is incident to
the above modular flag $F$.  The unique gallery $r_0$ from $c_0$ to
$-c_0$ incident to $F$ crosses the hyperplanes in this order:
\begin{align*}
&H^-_{11}, \\
&H^-_{12}, H^-_{22}, H^+_{12}, \\
&H^-_{23}, H^-_{13}, H^-_{33}, H^+_{13}, H^+_{23}, \\
&H^-_{34}, H^-_{24}, H^-_{14}, H^-_{44}, H^+_{14}, H^+_{24}, H^+_{34}, \ldots
\end{align*}
Choose Coxeter generators $S=\{s_0,s_1,\ldots,s_{n-1}\}$
for the hyperoctahedral group $W=B_n$ of signed permutations acting
on the coordinates of $\RR^n$, such that $s_i$ is the adjacent transposition $(i,i+1)$,
as before, and $s_0$ is the sign change in the first coordinate $x_1$.
Then this gallery $r_0$ corresponds to the following
reduced decomposition for $w_0$:
\begin{align*}
&(s_0, \\
&s_1, s_0, s_1, \\
&s_2, s_1, s_0, s_1, s_2, \\
&s_3, s_2, s_1, s_0, s_1, s_2, s_3, \ldots).
\end{align*}
\end{example}

\begin{remark}
Consider the following possible hypotheses on a real (central, essential)
hyperplane arrangement $\AAA$ and one of its chambers $c_0$:
\begin{enumerate}[(i)]
\item The chamber $c_0$ is a simplicial cone in the sense that
its walls have linearly independent normal vectors.
\item The weak order on the chambers $\CCC$, considered by 
Edelman~\cite{Edelman}, in which 
$c \leq c'$ when $L_1(c_0,c) \subseteq L_1(c_0,c')$, is a lattice.
\item $\AAA$ is \emph{simplicial}, meaning that
every chamber is a simplicial cone.
\item $\AAA$ is a real reflection arrangement.
\item The arrangement $\AAA$ is supersolvable, and the
chamber $c_0$ is incident to one of its modular flags.
\end{enumerate}

\noindent
Incorporating well-known results for
reflection arrangements with various results from
Bj\"orner, Edelman and Ziegler~\cite{BjornerEdelmanZiegler},
one has the following implications:
$$
\begin{array}{ccccccc}
(iv) & \rightarrow &  (iii) &              &       &             &      \\
     &             &      & \searrow&       &             &      \\
     &             &      &              & (ii) & \rightarrow & (i) \\
     &             &      & \nearrow&       &             &     \\
     &             & (v) &              &       &             &     
\end{array}
$$

Bearing in mind that Theorem~\ref{thm:supersolvable-theorem} assumes hypothesis (v), 
and the $d$-vertex-connectivity of the graph $G_2$ proven in 
\cite[Theorem 1.1]{AthanasiadisEdelmanReiner} assumes hypothesis (iii), 
it is reasonable to ask whether any of the extra hypotheses (i),(ii),(iii),(iv)
imply that the lower bound of $|L_2|$ for the diameter of the
graph $G_2$ is tight.
\end{remark}

\section{On the graphs $G(w)$}
\label{section:subgraphs}

In this section we show how the previous
methods generalize to the graph $G(w)$ of reduced
words for $w$ an element of a finite reflection group $W$, as discussed in the
Introduction.  Although these methods lead to some
bounds on the distance functions and diameters of $G(w)$ in
the groups of type $A, B$, we do not determine these diameters
exactly.

As in Section~\ref{section:first-two-graphs}, let $\AAA$ be a
(central, essential) hyperplane arrangement in $\RR^d$, with
set of chambers $\CCC$.  Given two chambers $c, c'\in \CCC$,
recall that $L_1(c,c')$ denotes the set of hyperplanes $H\in L_1 = \AAA$
that separate $c$ from $c'$, or equivalently, for which the chambers
$c/H, c'/H$ in the localized rank $1$ arrangement $\AAA_H$ are antipodal.
Further define
\begin{align*}
L_2(c,c')
:={}&\{ X \in L_2: c/X, c'/X \text{ are antipodal chambers in }\AAA_X \} \\
={}&\{ X \in L_2: \AAA_X \subseteq L_1(c,c') \}
\end{align*}
Denote by $\RRR(c,c')$ the set of all minimal galleries $r$ from
$c$ to $c'$.

Note that any minimal gallery $r$ in $\RRR(c,c')$ must cross each
of the hyperplanes in $L_1(c,c')$ exactly once, and is completely
determined by the linear order in which these hyperplanes are
crossed. For each codimension-two subspace $X$ in $L_2$, there are
two possibilities:
\begin{itemize}
\item
If $X \notin L_2(c,c')$ then the hyperplanes of $\AAA_X$
must be crossed by $r$ in a unique linear order, namely the
order in which one crosses the hyperplanes of $\AAA_X \cap L_1(c,c')$
when walking from $c/X$ to $c'/X$ in rank two.
\item
If $X \in L_2(c,c')$ then the hyperplanes of $\AAA_X$
can be crossed by $r$ in one of two possible linear orders,
as in Figure~\ref{fig:rank-two-figure}.
\end{itemize}

Consequently, given two minimal galleries $r, r'\in \RRR(c,c')$,
one can again speak of the separation set $L_2(r,r')$ as the
subset of codimension-two subspaces $X$ in $L_2(c,c')$ on which
$r, r'$ disagree\footnote{We apologize for the slightly confusing
overuse of ``$L_2$'' in our notations: we now have not only the set $L_2$, which is the
collection of all codimension-two intersection subspaces of
$\AAA$, but also two functions called $L_2(-,-)$, namely
$$
\begin{array}{rclcl}
\CCC &\times &\CCC &\longrightarrow &2^{L_2} \\
(c&,&c') &\longmapsto &L_2(c,c') \\
 & & & & \\
\RRR(c,c') &\times& \RRR(c,c') &\longrightarrow &2^{L_2(c,c')}\\
(r&,&r') &\longmapsto &L_2(r,r').
\end{array}
$$
We hope that context resolves
any confusion that arises within this section.}.

\begin{defn}[The graph $G_2(c,c')$]
\label{def:second-subgraph}
Given the arrangement $\AAA$ and two chambers $c, c'\in \CCC$,
define a graph $G_2(c,c')$ whose vertex set is
the set $\RRR(c,c')$ of minimal galleries $c$ to $c'$,
and having an edge between two galleries $\{r,r'\}$ exactly when $|L_2(r,r')|=1$.
\end{defn}

The previous discussion shows that the map
$L_2:\RRR(c,c') \times \RRR(c,c') \longrightarrow 2^{L_2(c,c')}$
provides a set-valued metric on $G_2(c,c')$, taking
values in $L_2(c,c')$.  That discussion also shows
that the generalization of the second assertion in
Proposition~\ref{prop:separation-sets-encode} still holds:  having fixed
a base gallery $r_0$ in $\RRR(c,c')$, then any other gallery
$r$ in $\RRR(c,c')$ is uniquely determined by its separation set
$L_2(r_0,r)$.  One has also this immediate consequence of
Propositions~\ref{prop:bipartiteness} and \ref{prop:set-lower-bound}.

\begin{corollary}
The graph $G:=G_2(c,c')$ is always bipartite, and
its distance function satisfies
$
d_G(r,r') \geq |L_2(r,r')|.
$
\end{corollary}

\begin{remark}
Whenever $c_0, c, c'\in \CCC$, satisfy $L_1(c_0,c) \subset
L_1(c_0,c')$, one can define an injection $\RRR(c_0,c)
\hookrightarrow \RRR(c_0,c')$: fix any minimal gallery $r$ from
$c$ to $c'$, and then concatenation with $r$ as a suffix
gives such an injection.  It is easily seen that this leads to an
embedding of the graph $G_2(c_0,c) \hookrightarrow G_2(c_0,c')$ as
a vertex-induced subgraph.

In particular, although the graphs $G_2(c_0,c)$ with their
set-valued metric do not in general have a $\ZZ_2$-equivariant involution,
they are always vertex-induced subgraphs of the graph $G_2(c_0,-c_0)$,
which does have such an involution.
\end{remark}

\begin{remark}
When $\AAA$ is the arrangement of reflecting hyperplanes for
a finite real reflection group $W$, the (simply) transitive
$W$-action on the chambers means that any graph $G_2(c,c')$
is isomorphic to the graph $G(w):=G_2(c_0,w(c_0))$ for some
fixed choice of a base chamber $c_0$ and some group element $w$
in $W$.  It is not hard to check that
this graph $G(w)$ is the graph of reduced words for $w$
discussed in the Introduction.

Furthermore, in this situation, the sets
\begin{align*}
L_1(w)&:=L_1(c_0,w(c_0)) \\
L_2(w)&:=L_2(c_0,w(c_0))
\end{align*}
have the following reflection group interpretations:
\begin{itemize}
\item
$L_1(w)$ is the usual (left-)inversion set of $w$, that
is, the collection of positive roots $\alpha_H$ for $W$ which are
sent to negative roots by $w^{-1}$.
\item
$L_2(w)$ is the collection of
rank two sub--root systems $\Phi_X$
having the property that $w^{-1}$ sends every positive root
in $\Phi_X$ to a negative root.
\end{itemize}
\end{remark}

One then has the following extension of
Theorem~\ref{thm:supersolvable-accessibility}.

\begin{theorem}
\label{thm:supersolvable-subgraph-accessibility}
Let $\AAA$ be a real (central, essential) hyperplane arrangement in $\RR^d$ which is
supersolvable, let $F:=\{X_i\}_{i=0}^d$ be a modular flag for $\AAA$,
let $c_0$ be any chamber incident to $F$, and let $c$ be any other chamber.

\begin{enumerate}[(i)]
\item
There is a unique minimal gallery
$r_0$ from $c_0$ to $c$ incident to $F$.
\item
This minimal gallery $r_0$ is an $L_2$-accessible vertex
for the graph $G_2(c_0,c)$ on the galleries $\RRR$ from $c_0$ to $c$.
\end{enumerate}
\end{theorem}

\begin{proof}
The proof of assertion (i) is by induction on $d$ exactly as in the proof
of Theorem~\ref{thm:supersolvable-accessibility}.  This unique
gallery $r_0$ is obtained by applying
induction to the $(d-1)$-dimensional supersolvable arrrangement $\AAA_\ell$,
lifting the unique gallery from $c_0/\ell$ to $c/\ell$ incident to $F/\ell$
in order to first cross all hyperplanes in
$L_1(c_0,c) \cap \AAA_\ell$.  One must complete it by then
crossing the hyperplanes in $L_1(c_0,c) \setminus \AAA_\ell$
in the linear order which is the restriction of the
one from Proposition~\ref{prop:modular-coatom-prop}(i).

The proof of assertion (ii) is also by induction on $d$ exactly as in the proof
of Theorem~\ref{thm:supersolvable-accessibility}.  One wishes to
show that for any gallery $r \neq r_0$ from $c_0$ to $c$, there
will be another gallery $r'$ having $L_2(r_0,r') \subset L_2(r_0,r)$
and $|L_2(r,r')| = 1$.
Again there are two cases, depending upon whether (Case 1) or not (Case 2)
the gallery $r$ crosses all the hyperplanes of $L_1(c_0,c) \cap \AAA_\ell$
before crossing any hyperplanes of $L_1(c_0,c) \setminus \AAA_\ell$.
There is no essential change in the proof of Case 1.  

In Case 2, one must note that the exhibited gallery $r'$ having
$L_2(r,r')=\{X\}$ satisfies in addition that $X$ lies in $L_2(c_0,c)$. 
This is immediate from the fact that $r$ crossed 
\emph{every} hyperplane of $\AAA_X$ on its way
from $c_0$ to $c$, so that $\AAA_X \subset L_1(c_0,c)$,
that is, $X$ lies in $L_2(c_0,c)$.
\end{proof}

\begin{example}
\label{type-A-example-3} Continuing
Example~\ref{type-A-example-2}, for the reflection arrangement $\AAA$ of
type $A_{5}$, choose as modular flag $F:=\{X_i\}_{i=0}^{5}$ where
$ X_i:=\{x_1=x_2=\cdots=x_{6-i} \}. $  The chamber $c_0:=\{x_{1} <
\cdots < x_{6}\}$ corresponding to the identity permutation is
incident to the above modular flag $F$. Let $w=316425\in \sym_6$.
The unique gallery $r_0$ from $c_0$ to $c:=w(c_0)$ incident to $F$ crosses the
hyperplanes in this order:
\begin{align*}
&H_{23}, H_{13},\\
&H_{24}, \\
&H_{56}, H_{26}, H_{46}. \\
\end{align*}
The corresponding sequence of chambers in the gallery is indexed
by the permutations
$$
123456 \overset{H_{23}}{-}
132456 \overset{H_{13}}{-}
312456 \overset{H_{24}}{-}
314256 \overset{H_{56}}{-}
314265 \overset{H_{26}}{-}
314625 \overset{H_{46}}{-}
316425, 
$$
which corresponds  to the following reduced decomposition for $w$
\begin{align*}
&(s_2, s_1,\\
&s_3, \\
&s_5, s_4, s_3). \\
\end{align*}

\end{example}

\begin{corollary}
\label{cor:supersolvable-subgraph-bounds}
Let $\AAA$ be a real (central, essential) hyperplane arrangement in $\RR^d$ which is
supersolvable, let $F:=\{X_i\}_{i=0}^d$ be a modular flag for $\AAA$,
and let $c_0$ be any chamber incident to $F$.
Let $c$ be any other chamber, and let $r_0$ be the unique
minimal gallery from $c_0$ to $c$ incident to $F$.

Then in the graph $G:=G_2(c_0,c)$ of minimal galleries from $c_0$ to $c$,
any two galleries $r, r'$ satisfy
$$
\begin{array}{rcl}
|L_2(r,r')| & \leq d_G(r,r') & \leq |L_2(r_0,r)| + |L_2(r_0,r')|  \\
            &                & \leq 2|L_2(c_0,c)| .
\end{array}
$$
Thus the diameter of $G$ is at most $2|L_2(c_0,c)|$.

In particular, when $\AAA$ is a reflection arrangement of type $A, B$ or
dihedral type $I_2(m)$, so that $G=G(w)$ for some $w\in W$,
any two reduced words $r, r'$ satisfy
$$
\begin{array}{rcl}
|L_2(r,r')| & \leq d_G(r,r') & \leq |L_2(r_0,r)| + |L_2(r_0,r')| \\
            &                & \leq 2|L_2(w)|.
\end{array}
$$
Thus the diameter of $G$ is at most $2|L_2(w)|$.
\end{corollary}

\begin{example}
Example~\ref{inaccessibility-example}
and the four words listed in \eqref{four-inaccessible-words}
show that both the upper and lower bounds
on $d_G(r,r')$ given in Corollary~\ref{cor:supersolvable-subgraph-bounds}
need not be tight, even for $G=G(w_0)$ in type $A_{n-1}$.

In light of Theorem~\ref{thm:supersolvable-theorem}
one might wonder whether the upper bound of $2|L_2(w)|$ for
the diameter of $G(w)$ in types $A, B$ can
be improved to an upper bound of $|L_2(w)|$.
However, we mention here some examples in types $A_3, B_3$
showing that even an upper bound of $|L_2(w)|$ for the diameter of $G(w)$
is not always tight.  

In type $A_3$, where $W$ is the symmetric group $\sym_4$, the permutation
$w=3412$ has two reduced words
\begin{align*}
\RRR(w)= \Bigl\{ & (s_2, s_3, s_1, s_2), \\
                    & (s_2, s_1, s_3, s_2) \Bigr\}
\end{align*}
so that the graph $G(w)$ is a single edge, having diameter $1$.  However,
in this case, $L_2(w)=\{ X_{13,24}, X_{14,23} \}$, so that $|L_2(w)|=2$.

One encounters a similar phenomenon in type $B_3$,
for the signed permutation $w$ that maps the standard
basis vectors $e_1,e_2,e_3$, respectively, to $-e_3,-e_2,-e_1$, respectively.
This $w$ has only two reduced words 
\begin{align*}
\RRR(w)= \left\{ \right. & (s_0, s_1, s_0, s_2, s_1 , s_0), \\
                    & \left. (s_0, s_1, s_2, s_0, s_1 , s_0) \right\}
\end{align*}
with respect to the Coxeter generators $S=\{s_0,s_1,s_2\}$ from
Example \ref{type-B-example-2}, 
so that the graph $G(w)$ is a single edge, having
diameter $1$.  But $L_2(w)$ consists of the three subspaces of
codimension-two of the form $\{ x_i=0, x_j + x_k = 0\}$ with $\{i,j,k\}=\{1,2,3\}$,
so that $|L_2(w)|=3$.
\end{example}

These examples suggest the following
conjecture\footnote{This conjecture has been checked only for the symmetric
groups $S_n$ with $n \leq 5$, and for the hyperoctohedral
groups $B_n$ with $n \leq 3$.}.

\begin{conjecture}\label{diameter-conjecture1}
For an element $w$ in the symmetric group $\sym_n$,
$$
\frac{1}{2}|L_2(w)|\le {\rm{diameter}}(G(w))\le |L_2(w)|.
$$
For an element $w$ in the hyperoctahedral group $B_n$,
$$
\frac{1}{3}|L_2(w)|\le {\rm{diameter}}(G(w))\le |L_2(w)|.
$$

\end{conjecture}

\begin{remark}
\label{Armstrong-remark}
The authors would like to thank an anonymous referee for pointing out the
following connection between 
Theorems~\ref{thm:supersolvable-accessibility}, 
\ref{thm:supersolvable-subgraph-accessibility}
and the work of Armstrong \cite{Armstrong} on {\it sorting orders}
in Coxeter groups.

Let $W$ be a real reflection group $W$, with reflection arrangement $\AAA$,
and assume $\AAA$ is supersolvable with a choice of modular flag $F$.
Theorem~\ref{thm:supersolvable-accessibility} 
shows there is a unique minimal gallery $r_0$ from $c_0$ to $-c_0$
incident to $F$, which corresponds to a particular reduced word ${\bf w}_0$
for the longest element $w_0$ in $W$.
Theorem~\ref{thm:supersolvable-subgraph-accessibility}
then shows that for any other chamber $c$, there is again a unique minimal gallery $r$
from $c_0$ to $c$ incident to $F$.
If $w$ is the unique element of $W$ for which $c=w(c_0)$,
then this minimal gallery $r$ corresponds to a particular reduced
word ${\bf w}$ for $w$.

Meanwhile in Armstrong's work, {\it any} choice of a reduced word ${\bf w}_0$
for $w_0$ induces, for each $w$ in $W$, a particular reduced word ${\bf w}$ for
$w$, which he calls the {\it ${\bf w}_0$-sorted word} for $w$.  Specifically,
${\bf w}$ is the lexicographically leftmost subword of ${\bf w}_0$
that gives a reduced word for $w$.

One can show that these two constructions are the same:
the word ${\bf w}$ corresponding to the gallery $r$ incident to $F$
is the same as the {\it ${\bf w}_0$-sorted word} for $w$.
This can be shown using an induction on the rank $d$, similar to the
one employed in the proof of Theorems~\ref{thm:supersolvable-accessibility} and
\ref{thm:supersolvable-subgraph-accessibility}, together
with the compatibility of parabolic coset
factorization with weak Bruhat orders, and 
\cite[Thm. 4.2]{Armstrong}.
\end{remark}

\section{Acknolwedgements}
The authors thank Anders Bj\"orner, Francesco Brenti, Patrick
Dehornoy, Axel Hultman, Arkadius Kalka, Jim Lawrence, Thomas McConville, 
John Sullivan, and G\"unter M. Ziegler for helpful discussions and references.  
They thank Rob Edman for computations in types $D_4, F_4$ 
checking the non-$L_2$-accessibility of lexicographcially first reduced words.
They thank Nathan Reading, Hugh Thomas, and an anonymous referee
for corrections of typos and errors in an earlier version.
Lastly, they thank another anonymous referee for helpful suggestions,
including the content of Remark~\ref{Armstrong-remark}.



\end{document}